\documentclass[11pt]{article}

\title{Scalar Representation and Conjugation of Set--Valued Functions}

\author{Carola Schrage
\footnote{\href{mailto:carola.schrage@mathematik.uni-halle.de}{carola.schrage@mathematik.uni-halle.de}} }
\date{{\small \today}}

\usepackage{amsmath, amssymb, enumerate}

\newtheorem{theorem}{Theorem}
\newtheorem{corollary}[theorem]{Corollary}
\newtheorem{remark}[theorem]{Remark}
\newtheorem{lemma}[theorem]{Lemma}
\newtheorem{definition}[theorem]{Definition}
\newtheorem{proposition}[theorem]{Proposition}
\newtheorem{example}[theorem]{Example}

\numberwithin{equation}{section}  
\numberwithin{figure}{section}    
\numberwithin{table}{section}     
\numberwithin{theorem}{section}

\newcommand{\proof}{ {\sc Proof.\quad}}
\newcommand{\pend}{ \hfill $\square$ \\}

\newcommand{\of}[1]{\ensuremath{\left( #1 \right)}}
\newcommand{\cb}[1]{\ensuremath{ \left\{ #1 \right\} }}

\newcommand{\st}{\,|\; }

\newcommand{\vp}{\ensuremath{\varphi}}

\renewcommand{\P}{\ensuremath{\mathcal{P}}}

\newcommand{\G}{\ensuremath{\mathcal{G}}}

\newcommand{\lel}{\preccurlyeq}
\newcommand{\leu}{\curlyeqprec}

\newcommand{\R}{\mathrm{I\negthinspace R}}
\newcommand{\OLR}{\overline{\mathrm{I\negthinspace R}}}
\newcommand{\N}{\mathrm{I\negthinspace N}}

\newcommand{\dom}{{\rm dom \,}}
\newcommand{\epi}{{\rm epi \,}}

\newcommand{\gr}{{\rm graph \,}}
\newcommand{\cl}{{\rm cl \,}}
\newcommand{\co}{{\rm co \,}}
\newcommand{\cone}{{\rm cone\,}}
\newcommand{\Int}{{\rm int\,}}

\newcommand{\dual}{{\rm c \,}}

\newcommand{\isum}{{+^{\negmedspace\centerdot\,}}}
\newcommand{\ssum}{{+_{\negmedspace\centerdot\,}}}
\newcommand{\idif}{{-^{\negmedspace\centerdot\,}}}
\newcommand{\sdif}{{-_{\negmedspace\centerdot\,}}}
\newcommand{\isquare}{{\square^{\centerdot\,}}}
\newcommand{\ssquare}{{\square_{\centerdot\,}}}

\newcommand{\triup}{{\rm \vartriangle}}
\newcommand{\trido}{{\rm \triangledown}}

\usepackage[
colorlinks=true, linkcolor=blue, citecolor=blue, urlcolor=blue, filecolor=blue
]{hyperref}

\begin{document}
\maketitle
\begin{abstract}
We are considering functions with values in the power set of a pre--ordered, separated locally convex space with closed convex images. To each such function, a family of scalarizations is given which completely characterizes the original function.  A concept of a Legendre--Fenchel conjugate for  set--valued functions is introduced and identified with the conjugates of the scalarizations. To the set--valued conjugate, a full calculus is provided, including a biconjugation theorem, a chain rule and weak and strong duality results of Fenchel-Rockafellar type.


{\bf Classcode:} 49N15, 54C60, 90C46
\end{abstract}

\section{Introduction}

In this paper, we will introduce a new duality theory for set--valued functions. We will apply our theory proving a biconjugation theorem, a sum-- and a chain--rule and weak and strong duality results of Fenchel--Rockafellar type.
Apart from the purely academic interest, investigating set--valued functions  is motivated by their applicability in  financial mathematics, compare \cite{CasMo06, HamelHeyHoe07, JouMeTou04}, and in vector optimization.

The common lack of infima and suprema in vector spaces makes it obvious that 'vectorial' constructions are in general not appropriate to solve minimization problems on a vector space, compare \cite{Hamel09}.
Even if investigations are restricted to order complete, partially  ordered spaces or subsets, least upper or greatest lower bounds can be 'far away' and have little in common with the function in question.

An approach to solve this dilemma is to search for a set of minimal elements and hereby transform the original vector--valued problem into a set--valued problem.

Apparently, as pointed out by J.Jahn, 'the best soccer team is not necessarily the team with the best player', thus we understand an optimal solution of a set--valued minimization problem to be a set rather than a single point. Consequently, we investigate set--valued functions understanding the images to be elements of the power set of a vector space.

The basic idea is to extend the order on the space $Z$, given by a convex cone $C\subseteq Z$ to an order on $\P(Z)$ in an optimal way, compare \cite{Kuroiwa97, HamelHabil05, Jahn04} and identify a function $f:X\to\P(Z)$ with its epigraphical extension $f_C(x)=f(x)+C$.
This extension can of course be done for vector--valued functions, thus as a special amenity our theory includes vector--valued functions as a special class of set--valued functions.
We identify the subset $\P^\triup\subseteq \P(Z)$ as the set of all elements $A=A+C\in\P(Z)$. This set is an order complete lattice w.r.t. $\supseteq$ and contains all images of functions of epigraphical type. Equipped with an appropriate addition and multiplication with nonnegative reals, $\P^\triup$ is a so--called $\inf$--residuated conlinear space, compare \cite{HamelHabil05, HamelSchrage10}, thus supplies sufficient algebraic structure to introduce coherent definitions of conjugates, directional derivatives and subdifferentials. At present, we will concentrate on the investigation of conjugates, more can be found in \cite{Diss}, compare also \cite{demyanov1986approximation} and especially \cite{Loehne11Buch} for closely related approaches.

The introduced notion of conjugate functions turns out to be an example of a $(*,s)$--duality, see \cite{GetanMaLeSi}, and many properties from the well--known scalar conjugate can be rediscovered in our setting.

In \cite{Hamel09, Hamel10}, duality results are proven directly, using separation arguments in the space $X\times Z$. In contrast, we introduce a family of scalar functions related to and completely characterizing a set--valued function. The set--valued conjugate is also determined by the family of conjugates of the scalarizations of the original function.

Because of the one-to-one correspondence between set--valued functions and their conjugates on the one hand, the family of scalarizations and their conjugates on the other, we can utilize results from the scalar theory to obtain results for the set--valued case.

The article is organized as follows. 
In Section \ref{sec:ConlinSpace} we resume basic definitions and facts about conlinear spaces, residuation, extended real--valued and set--valued functions. Apart from an extended chain--rule, everything presented in this section can be found in the literature, compare e.g. \cite{galatos2007residuated, HamelHabil05, HamelSchrage10} and also \cite{IoffeTikhomirov79, Zalinescu02}.
In Section \ref{sec:Scalarization}, a family of scalar functions associated to a set--valued function is identified.
This family of scalarizations proves to characterize its assigned set--valued counterpart and will prove useful in later course to provide a concise approach to a theory of convex analysis for set--valued functions.
The definitions of conjugate and biconjugate are provided in Section \ref{sec:Conjugation}, while a selection of duality results is presented in Section \ref{sec:Duality}.

\section{Preliminaries; Residuation and Conlinear Structure}\label{sec:ConlinSpace}

\subsection{General conlinear spaces}\label{subsec:GeneralConlinStructure}

In \cite[Section 2.1.2]{HamelHabil05}, the concept of a conlinear space has been introduced. A set $Z=\of{Z,+,\cdot}$ is called a conlinear space, if $(Z,+)$ is a commutative monoid with neutral element and for all $z,z_1,z_2\in Z$, $r,s\in \R_+$ holds $r(z_1+z_2)=rz_1+rz_2$, $r(sz)=(rs)z$ and $1z=z$, $0z=0$.

If $(Z,\leq)$ is a complete lattice and $\leq$ is compatible with the algebraic structure in $\of{Z,+,\cdot}$, then $(Z,+,\cdot,\leq)$ is called an order complete conlinear space.

If additionally the operation $+:Z\times Z\to Z$ satisfies $(x+\inf M)=\inf\limits_{m\in M} \of{x+m}$, then the $\inf$--residuation or $\inf$--difference  $\idif:Z\times Z\to Z$ defined by
\begin{equation}\label{eq:InfResBasic}
\forall z_1, z_2\in Z:\quad z_1\idif z_2 := \inf\cb{t\in Z\st z_1\leq z_2+t}
\end{equation}
replaces the usual difference operation. It holds $z_1\leq z_2+ (z_1\idif z_2)$ and $Z=\of{Z,+,\cdot,\leq}$ is called an $\inf$--residuated order complete conlinear space.

Dually, if $(x+\sup M)=\sup\limits_{m\in M} \of{x+m}$, then the $\sup$--residuation ($\sup$--difference) $\sdif:Z\times Z\to Z$ defined by
\[
\forall z_1, z_2\in Z:\quad z_1\sdif z_2 = \sup\cb{t\in Z\st z_2+t\leq z_1}
\]
replaces the difference operation. It holds
$z_2+ (z_1\sdif z_2)\leq z_2$ and $Z=\of{Z,+,\cdot,\leq}$ is called a $\sup$--residuated order complete conlinear space.

A conlinear $\inf$--residuated order complete space as image space supplies sufficient structure to apply the concepts introduced in \cite{GetanMaLeSi} and to define convexity of functions, as multiplication with positive reals is defined as was done in \cite{Hamel09}.
Also, the structure carries over from a space $Z$ to the set of all subsets, $\P(Z)$, provided the algebraic and order relations are extended to relations on $\P(Z)$ in an appropriate way, see \cite[Theorem 13]{HamelHabil05} and \cite{HamelSchrage10} for more details.

Residuation seems to be introduced by Dedekind, \cite[p. 71]{Dedekind1863},\cite[p. 329-330]{Dedekind1872}, compare also \cite[XIV, \S 5]{Birkhoff40}, \cite[chap. XII]{Fuchs66} or \cite[Chap. 3]{galatos2007residuated} on $\sup$--residuation. In \cite[Lemma 3.3]{galatos2007residuated} the concept of $\inf$--residuation is indicated.
The corresponding residuations may substitute the difference operation on $\R$.

\subsection{The extended real numbers}\label{subsec:ExtReals}

The set $\R$ is extended to $\OLR$ in the usual way by adding two elements $+\infty$ and $-\infty$. Setting
\begin{align*}
  & \forall t\in\R:\quad -\infty \leq t \leq +\infty;\\
  & \inf\emptyset=+\infty,\quad \sup\emptyset=-\infty,
\end{align*}
the set $(\OLR,\leq)$ is an order complete lattice.
Multiplication with nonnegative reals is extended to $\OLR$ by setting $t\cdot(\pm\infty)=\pm\infty$, if $t>0$ and $0\cdot(\pm\infty)=0$.
The addition $+:\R\times \R\to \R$ on $\R$ admits two distinct extensions to an operator on $\OLR$, the $\inf$-- and $\sup$--addition.
\begin{align}\label{eq:InfSupAddition}
 \forall r,s\in\OLR:\quad  r\isum s & = \inf\cb{a+b\st a,b\in \R,\, r\leq a,\, s\leq b};\\
 \label{eq:InfSupAdditionII}
                            r\ssum s & = \sup\cb{a+b\st a,b\in \R,\, a\leq r,\, b\leq s}.
\end{align}
These constructions have been introduced in \cite[p. 329-330]{Dedekind1872} in order to extend the addition from the set of rational to the set of real numbers.

Both $\inf$-- and $\sup$--addition are commutative and compatible with the usual order $\leq$ on $\OLR$.
The greatest element $+\infty$ dominates the $\inf$--addition $\isum$, $-\infty$ dominates the $\sup$--addition,
\begin{align}\label{eq:inf-sup-Addition}
\forall r\in\OLR:\quad (+\infty)\isum r=+\infty,\; (-\infty)\ssum r=-\infty.
\end{align}
The sets
\begin{align*}
  \R^\triup  =(\OLR,\isum,\cdot,\leq),\quad \R^\trido  =(\OLR,\ssum,\cdot,\leq)
\end{align*}
are $\inf$-- and $\sup$--residuated complete conlinear spaces, respectively.
The suffix $\triup$ indicates $\inf$--residuation and  dually the suffix $\trido$ indicates $\sup$--residuation.
Multiplication with $-1$ is defined as $(-1)\cdot(\pm\infty)=\mp\infty$. We abbreviate $(-1)r=-r$, when no confusion can arise.
Obviously, for $M\subseteq\OLR$ and $r,s\in\OLR$ the following is satisfied.
\begin{align*}
& (-1)\inf M =\sup(-1) M; \\
& r\isum (-s)= r\sdif s,\quad -(r\isum s)= (-s)\ssum (-r);\\
& r\ssum (-s)=r\idif s,\quad -(r\ssum s)=(-s)\isum (-r).
\end{align*}
Multiplication with $-1$ is a duality between $\R^\triup$ and $\R^\trido$.

\subsection{Extended real--valued functions}\label{subsec:ExtendedValuedFunctions}

Let $X$ and $Y$ be locally convex separated spaces with topological duals $X^*$ and $Y^*$ and $g:X\to\R^\triup$ a function. Multiplication with $-1$ transfers $g$ to $-g:X\to\R^\trido$. We will concentrate on the first type of functions in the sequel, keeping in mind that for the second class symmetric results can be proven.
To $g:X\to\R^\triup$, the epigraph and effective domain of $g$ are defined as usual
\begin{align*}
\epi g  = \cb{(x,r)\in X\times \R\st g(x)\leq r},\quad
\dom g = \cb{x\in X\st g(x)\neq +\infty}.
\end{align*}

The set $(\R^\triup)^X = \cb{g:X\to\R^\triup}$ equipped with the point--wise addition, multiplication with nonnegative reals and order relation is an order complete, $\inf$--residuated conlinear space.
The $\inf$--convolution with respect to either $\inf$--addition or $\sup$--addition is denoted as $\isquare$ and $\ssquare$, respectively.
We denote
\begin{align}\label{eq:fT_and_Tg2a}
 \forall x\in X:\quad (fT)(x) =f(Tx); \quad (Tg)(y)  =\inf\limits_{Tx=y}g(x)
\end{align}
for a linear continuous operator $T:X\to Y$,  $f:Y\to \R^\triup$ and $g:X\to \R^\triup$.

A function $g:X\to\R^\triup$ is said to be closed (convex), if its epigraph is a closed (convex) set.
It is subadditive, if $\epi g$ is closed under addition, i.e. $\epi g + \epi g\subseteq \epi g$  and positively homogeneous, if $\epi g$ is a cone, i.e.
\[
\epi g = \cone\of{\epi g}=\cb{t(x,r)\in X\times\R\st t>0, (x,r)\in \epi g}.
\]
A positively homogeneous convex function is called sublinear. Equivalently, a function is sublinear if and only if it is positively homogeneous and subadditive.

A function $g:X\to\R^\triup$ is proper if $\dom g\neq\emptyset$ and for all $x\in X$ exists a $t\in \R$ such that $(x,t)\notin \epi g$.

A well known separation theorem is the following corollary of the Hahn Banach Theorem.

\begin{lemma}\cite[Corollary 1.4]{EkelandTemam}\label{lem:Separation}
In a locally convex separated space $Z$, every closed convex set is the intersection of all closed (affine) half--spaces containing it.
\end{lemma}

As each such affine half--space can be expressed via a non--zero affine continuous function, $H_\alpha(z^*)=\cb{z\in Z\st \alpha\leq -z^*(z)}$, it is immediate that a function $g:X\to\OLR$ is closed, convex and either proper or constant $+\infty$ or $-\infty$, if and only if it is the supremum of its affine minorants, in which case the epigraph of $g$ as a subset of $X\times\R$ is the intersection of all $H_\alpha(x^*,-1)$, where $x^*(x)+\alpha\leq g(x)$ for all $x\in X$. A similar result will be provided for set--valued functions in Theorem \ref{thm:DualRepresentation}.
Equivalently to Lemma \ref{lem:Separation}, one can state that in a locally convex separated space $Z$, whenever there is a closed convex set $M\subseteq Z$ and an element $z\notin M$, then it exists a non--zero element $z^*\in Z^*\setminus\cb{0}$ and $t\in \R$, such that 
\begin{align}\label{eq:separation}
-z^*(z)<t\leq\inf\limits_{m\in M}-z^*(m).
\end{align}
Equation \eqref{eq:separation} will prove useful in various occasions throughout this paper.


The conjugate and biconjugate of $g:X\to\R^\triup$ are given by
\begin{align}\label{eq:ScalarConjugate}
    \forall x^*\in X^*:\quad g^*(x^*)& =\sup\limits_{x\in X}\of{x^*(x)\idif g(x)};\\
\label{eq:ScalarBiconjugate}
    \forall x\in X:\quad  g^{**}(x) & =\sup\limits_{x^*\in X^*}\of{x^*(x)\idif g^*(x^*)}.
\end{align}

For completeness we cite the fundamental duality formula, Theorem \ref{thm:ScalarFundamentalDuality}
and sum up these considerations with an extended chain--rule for scalar valued functions.
To our ends, it is of no consequence that in \cite{Zalinescu02} the product of $0\cdot(+\infty)$ is defined differently from our definition.

\begin{theorem}\cite[Theorem 2.7.1(iii)]{Zalinescu02}\label{thm:ScalarFundamentalDuality}
  Let $X$ and $Y$ be topological linear spaces with topological duals $X^*$ and $Y^*$ and $h:X\times Y\to\OLR$ a proper convex function with $(x_0,0)\in\dom h$ for some $x_0\in X$. If $h(x_0,\cdot):Y\to\OLR$ is continuous in $0\in Y$, then
  \begin{align*}
    \inf\limits_{x\in X} h(x,0) & =\sup\limits_{y^*\in Y^*} 0\idif h^*(0,y^*);\\
    \exists y^*\in Y^*:\quad \inf\limits_{x\in X} h(x,0) & =0\idif h^*(0,y^*).
  \end{align*}
  Moreover, if $\bar x\in X$, then
  \begin{align*}
   & \exists \bar y^*\in Y^*:\quad (0,\bar y^*)(\bar x, 0)\idif h(\bar x,0)= h^*(0,\bar y^*)\\
   \Leftrightarrow\quad & h(\bar x,0)=  \inf\limits_{x\in X} h(x,0).
  \end{align*}
\end{theorem}

\begin{theorem}[Scalar Chain Rule ]\label{thm:ScalarChainRule}
  Let $Y$ be another  topological linear space with topological dual $Y^*$, $g:X\to\R^\triup$, $f:Y\to\R^\triup$ two functions and $T:X\to Y$, $S:Y\to X$ linear continuous operators.
\begin{enumerate}[(a)]
  \item
   The conjugate of $\of{g\isquare Sf}$ is the function
\[
\forall x^*\in X^*:\quad \of{g\isquare Sf}^*(x^*)=\of{ g^*\ssum f^*S^*}\of{x^*}
                        \leq \of{ g^*\isum f^*S^*}\of{x^*} .
\]

\item
   The conjugate of $\of{g\isum fT}$ is dominated by $\of{g^*\ssquare T^*f^*}:X^*\to\R^\triup$,
\[
\forall x^*\in X^*:\quad
\of{g\isum fT}^*(x^*)\leq \of{g^*\ssquare T^*f^*}\of{x^*} \leq \of{g^*\isquare T^*f^*}\of{ x^*}.
\]

\item
    If additionally $g$ or $f$ is the constant mapping $+\infty$, then
\begin{align*}
  \forall x^*\in X^*:\quad   -\infty=\of{g\isum fT}^*(x^*)
            & = \of{g^*\ssquare T^*f^*}\of{x^*};\\
 \forall x^*\in X^*,\;\forall y^*\in Y^*:\quad
         \of{g\isum fT}^*(x^*)
            & =  g^*\of{x^*-T^*y^*}\ssum f^*\of{y^*}.
\end{align*}

\item
    If  $(fT)(x_0)=-\infty$ is satisfied for some $x_0\in \dom g$ or if both $f$ and $g$ are proper, convex and  $f$ is continuous in  a point in  $T(\dom g)$, then
\begin{eqnarray*}
  \forall x^*\in X^*:\quad -\infty\neq &\of{g\isum fT}^*(x^*)&= \of{g^*\isquare T^*f^*}\of{ x^*};\\
  \exists y^*\in Y^*:\quad &\of{g\isum fT}^*(x^*)& =  g^*\of{x^*-T^*y^*}\isum f^*\of{y^*}.
\end{eqnarray*}
\end{enumerate}
\end{theorem}
\proof
\begin{enumerate}[(a)]
\item
    We apply calculus rules in the residuated spaces $\R^\triup$ and $\R^\trido$ to obtain the following for all $x^*\in X^*$.
\begin{align*}
  \of{g\isquare Sf}^*\of{x^*} & =
        \sup\limits_{x\in X} \of{ x^*(x)\idif \inf\limits_{\bar x\in X}\of{ g(x-\bar x)\isum
                                                                \inf\limits_{Sy=\bar x}f(y) } }\\
    & = \sup\limits_{x\in X,\, y\in Y}
            \of{ \of{x^*(x-Sy)\idif g(x-Sy) } \ssum \of{ S^*x^*(y)\idif f(y) } }\\
    & = g^*(x^*)\ssum (f^*S^*)(x^*)
\end{align*}
and the second inequality holds by definition of $\ssum$ and $\isum$.

\item
   First we consider $T^*y^*=x^*$:
\[
(fT)^*(x^*)=\sup\limits_{x\in X}\of{y^*(Tx)\idif f\of{T(x)} }
           \leq \sup\limits_{y\in Y}\of{y^*(y)\idif f(y)} = f^*(y^*).
\]
We apply results from residuation theory to prove the following for all $x^*, \bar x^*\in X^*$.
\begin{align*}
  \of{g\isum fT}^*(x^*) & = \sup\limits_{x\in X}\of{
        \of{x^*(x)+ \bar x^*(x)-\bar x^*(x) }\ssum (-1)\of{g\of{x}\isum (fT)\of{x}} }\\
            & = \sup\limits_{x\in X}\of{
        \of{(x^*-\bar x^*)(x)\idif g\of{x}}\ssum \of{ \bar x^*(x) \idif(fT)\of{x}} }\\
            & \leq g^*(x^*-\bar x^*) \ssum (fT)^*(\bar x^*)\\
            & \leq  g^*(x^*-\bar x^*) \ssum T^*f^*(\bar x^*)
\end{align*}
and the second inequality holds by definition of $\ssum$ and $\isum$.

\item
    and the improper case of $(d)$ are proven by direct calculation, while the proper case is classic, compare e.g. \cite[Chap. 3, \S 3.4 Theorem 1]{IoffeTikhomirov79}.
\pend
\end{enumerate}

\begin{example}
  Under the assumptions of Theorem \ref{thm:ScalarChainRule}, let $g\equiv +\infty$ and $f(Tx_0)=-\infty$ with $x_0\in X$. Then
\begin{align*}
  \forall x\in X:&\quad \of{g\isquare Sf}(x) =  \of{g\isum fT}(x) = +\infty;\\
  \forall x^*\in X^*:&\quad g^*(x^*)=-\infty; \quad f^*S^*(x^*) = T^*f^*(x^*) =+\infty. 
 \end{align*}
Therefore
\begin{align*}
-\infty&=\of{g\isquare Sf}^*\of{x^*}  \leq \of{ g^*\isum f^*S^*}\of{x^*}=+\infty;\\
-\infty&=\of{g\isum fT}^*(x^*)\leq \of{g^*\isquare T^*f^*}\of{ x^*}=+\infty
\end{align*}
for all $x^*\in X^*$. The statements of Theorem \ref{thm:ScalarChainRule}(a) and (c) do not apply with equality for the $\inf$--addition on the right hand side.
In general, $\sup$--addition is dominated by $\inf$--addition and both operators coincide if neither addend is $-\infty$.
\end{example}

\subsection{The set $\P^\triup$}\label{sec:Set-Valued}

In the sequel, we will consider $X$, $Y$ and $Z$ to be locally convex separated  spaces with topological duals $X^*$, $Y^*$, $Z^*$ and investigate on functions from $X$ to a subset $\P^\triup$ of the power--set  $\P\of{Z}$ which is an order complete, $\inf$--residuated conlinear space. The suffix $\triup$ indicates $\inf$--residuation and  dually the suffix $\trido$ indicates $\sup$--residuation.

On $Z$, a reflexive and transitive order is given by a convex cone $C\subseteq Z$ with $\cb{0}\subsetneq C$, setting $z_1\leq z_2$, iff $z_2-z_1\in C$. The negative dual cone of $C$ is $C^-$, defined by
\[
C^- = \cb{z^*\in Z^*\st  \forall c\in C:\; z^*(c)\leq 0}
\]
and we assume $\cb{0}\subsetneq C^-$. To an element $z^*\in Z^*$ we define $H(z^*)=\cb{z\in Z\st  0\leq -z^*(z)}$.
Obviously, $C\subseteq H(z^*)$, if $z^*\in C^-$ and $\cl C= \bigcap\limits_{z^*\in C^-}H(z^*)$.
We do not assume  the topological interior of $C$ to be nonempty.

On $\P(Z)$, we introduce an algebraic structure by defining
\begin{align*}
  \forall A, B\in \P(Z):\quad A+ B & =\cb{a+b\st  a\in A,\, b\in B},\\
  \forall t\in\R\setminus\cb{0}:\quad tA & =\cb{ta\st  a\in A},
\end{align*}
the Minkowsky sum of two sets and the product of a set with a real number $t\neq 0$.
By convention, $A+\emptyset=\emptyset+A=\emptyset$ and $t\emptyset=\emptyset$ for all $A\in\P(Z)$ and $t\neq 0$ while $0 A=\cb{0}$ for all $A\in\P(Z)$.
We abbreviate $z+A=\cb{z}+A$ and $-A=(-1)A$ as well as $A-B=A+(-1)B$ for $A, B\in \P(Z)$ and $z\in Z$.

The order relation in $Z$ can be extended in two distinct ways by setting
\begin{align*}
  \forall A, B\in\P(Z):\quad & A\lel_C B\quad \Leftrightarrow B\subseteq A+C;\\
  \forall A, B\in\P(Z):\quad & A\leu_C B\quad \Leftrightarrow A\subseteq B-C,
\end{align*}
see \cite{Kuroiwa97, HamelHabil05, LoehneDiss}.

Two sets $A, B\in\P(Z)$ are equivalent with respect to $\lel_C$, iff  $A+C=B+C$. We identify the set $\P^\triup(Z,C)=\cb{A\in \P(Z)\st  A=A+C}$ and abbreviate $\P^\triup=\P^\triup(Z,C)$, if no confusion can arise. In $\P^\triup$, $A\lel_C B$ is equivalent to $B\subseteq A$. The set $\P^\triup$ is a complete lattice, infimum and supremum of a nonempty set $\mathcal M\subseteq \P^\triup$ are given by
\begin{align*}
  \inf \mathcal M  = \bigcup\limits_{M\in \mathcal M} M \in \P^\triup;\quad
  \sup \mathcal  M  = \bigcap\limits_{M\in\mathcal M} M \in \P^\triup
\end{align*}
and by convention $\inf\emptyset= \emptyset$ and $\sup\emptyset=Z$. The greatest element in $\P^\triup$ is $\emptyset$, the smallest $Z$.
For any $\mathcal M\subseteq \P^\triup$ and $A\in \P^\triup$,  
\begin{align*}
\inf(A+\mathcal M)  =      A+\inf \mathcal M;\quad
\sup (A+\mathcal M)  \supseteq A+\sup \mathcal M
\end{align*}
is satisfied.
Altering the multiplication with $0$ to $0A=C$ for all $A\in \P^\triup$, the set $\P^\triup$ together with the Minkowsky sum, the altered multiplication with nonnegative reals and the order relation $\supseteq$ is an order complete, $\inf$--residuated conlinear space. According to Equation \eqref{eq:InfResBasic}, the corresponding $\inf$--residuation is given by
\[
A\idif B= \inf\cb{M\in \P^\triup\st  A\lel_C B+ M}=\cb{z\in Z\st  A\supseteq B+z},
\]
In case the ordering cone is set to $C=\cb{0}$, this operation is sometimes called the star difference, compare \cite[Section 4]{penot2011} and the references therein.

\begin{lemma}\label{lem:Review}
For $z^*\in Z^*$, $\alpha\in\OLR$ and $A\subseteq Z$, set $\sigma(z^*|A)=\sup\cb{z^*(a)\st  a\in A}$,  $\sigma(z^*|\emptyset)=-\infty$ and
\begin{equation*}
H_{\alpha}(z^*)=\cb{z\in Z\st \alpha \leq -z^*(z)}.
\end{equation*}
Then $H(z^*)=H_0(z^*)$, $H_{\alpha\isum(-\sigma(z^*|A))}(z^*)= \cl(A+H_\alpha(z^*))$ and $H_{\alpha}(z^*)\idif A = H_{\alpha\idif(-\sigma(z^*|A))}(z^*)$.
\end{lemma}
\proof
It holds 
\[
\cl\of{A+H_\alpha(z^*)}=\cl\bigcup\limits_{a\in A}\of{a+H_\alpha(z^*)}
\]
and either $A=\emptyset$ and $H_{\alpha\isum(-\sigma(z^*|A))}(z^*)= \cl(A+H_\alpha(z^*))=\emptyset$ , or $\of{a+H_\alpha(z^*)}=H_{\alpha\isum(-z^*(a))}(z^*)$ is satisfied for all $a\in A$. Moreover,
\[
\cl\bigcup\limits_{a\in A}H_{\alpha\isum(-z^*(a))}(z^*)=H_{\inf\limits_{a\in A}\of{\alpha\isum(-z^*(a))}}(z^*)
\]
is true and by the calculus rules in $\R^\triup$, 
\[
\inf\limits_{a\in A}\of{\alpha\isum(-z^*(a))}= \alpha\isum(-\sigma(z^*|A)),
\]
hence $H_{\alpha\isum(-\sigma(z^*|A))}(z^*)= \cl(A+H_\alpha(z^*))$.
Likewise, by applying the residuation rules one proves
\begin{equation*}
H_{\alpha}(z^*)\idif A = H_{\alpha\idif(-\sigma(z^*|A))}(z^*).
\end{equation*}
\pend

We abbreviate
\begin{align}
  \P^\triup=\of{ \P^\triup,+,\cdot,\supseteq};\quad \P^\trido=\of{\cb{A-C\st  A\in\P(Z)}, +, \cdot,\subseteq}.
\end{align}
Multiplication with $-1$ is a duality between $\P^\triup$ and the order complete, $\sup$--residuated conlinear space $\P^\trido$.
With convex duality in mind, the space $\P^\triup$ is more appropriate, a closer study of both extensions can be found in \cite{HamelHabil05, loehne2005optimization}.

\subsection{Set--valued functions}\label{sec:SetValuedFunctios}

\begin{definition}\label{def:EpiHypoDomSetVal}
The graph of a function $g:X\to\P(Z)$ is defined as
\begin{align*}
  \gr g   & = \cb{(x,z)\in X\times Z\st  z\in g(x)}.
\end{align*}
If $g:X\to \P^\triup$, the domain and epigraph of $g$ are given by
\begin{align*}
  \dom g   = \cb{x\in X\st  g(x)\neq\emptyset},\quad
  \epi g   = \cb{(x,z)\in X\times Z\st  z\in g(x)+C}.
\end{align*}
\end{definition}

If $g:X\to \P^\triup$, then $\gr g=\epi g$, motivating the notion \emph{epigraphical type function} for functions $g:X\to \P^\triup$.

The set $(\P^\triup)^X=\cb{g:X\to \P^\triup}$  equipped with the point--wise addition, multiplication, order relation and $\inf$--residuation is an $\inf$--residuated order complete conlinear space.
We denote the $\inf$--convolution of $f,g:X\to \P^\triup$ by
\begin{align}\label{eq:InfConvolution}
\forall x\in X:\quad (f\isquare g)(x)= \inf\limits_{y\in X}\of{f(x-y)+ g(y)}
\end{align}
and 
\begin{align}\label{eq:fT_and_Tg2}
 \forall x\in X:\quad (fT)(x) =f(Tx); \quad (Tg)(y)  =\inf\limits_{Tx=y}g(x)
\end{align}
for a linear continuous operator $T:X\to Y$, $f:Y\to \P^\triup$ and $g:X\to \P^\triup$. The definitions in \eqref{eq:InfConvolution} and \eqref{eq:fT_and_Tg2} can be found in \cite{Hamel09,Diss}.

\begin{definition}\label{def:ConvexProperties}
  A function $g:X\to \P^\triup$ is called positively homogeneous, iff $\epi g$ is a cone,
  convex or closed, iff $\epi g$ is convex or closed, subadditive iff $\epi g$ is closed under addition and sublinear, iff $\epi g$ is a convex cone.
\end{definition}
Again, a positively homogeneous function is sublinear, if and only if it is subadditive.

  If a function $g:X\to \P^\triup$ is convex or closed, then especially for each $x\in X$ the set $g(x)$ is convex or closed.
  Defining $(\cl \co g):X\to \P^\triup$ by setting $\epi (\cl \co g)=\cl \co \of{\epi g}$, the function $(\cl \co g)$ maps into the set $\G^\triup$ of all convex, closed sets $A\in \P(Z)$ with $A=\cl\co(A+C)$.
  Altering the addition of sets to $A\oplus B=\cl(A+B)$ and multiplication with $0$ to $0A=\cl C$, then $\G^\triup$ with the altered addition, multiplication and the order relation $\supseteq$ is an $\inf$--residuated order complete conlinear space, the infimum in $\G^\triup$ of a set $\mathcal A\subseteq \G^\triup$ is $\inf_{\G^\triup}\mathcal A=\cl\co\bigcup\limits_{A\in \mathcal A}A$. 
This space has been used as image space in \cite{Hamel09,Diss}.

\begin{definition}\label{def:PropernesForSVFunctions}
  A function $g:X\to \P^\triup$ is called proper, iff $\dom g\neq \emptyset$ and there is no $x\in X$ with $g(x)=Z$.
  A function $g$ is called $z^*$--proper with $z^*\in Z^*$, iff
$x\mapsto-\sigma(z^*|g(x))$ is a proper function.
\end{definition}
Each function with nonempty domain is $0$-proper and if $g:X\to\P^\triup$ is $z^*$--proper, then $z^*\in C^-$.
If a closed convex function $g:X\to \P^\triup$ is improper, then $g(x)=Z$ holds for all $x\in\dom g$ and $\dom g$ is a closed convex set in $X$. Likewise, if a closed convex function is $z^*$--improper, then $(g(x)-H(z^*))\setminus(g(x)+H(z^*))=\emptyset$ is satisfied for all $x\in\dom g$ and $\dom g$ is closed and convex, compare \cite[Proposition 5]{Hamel09}.

\section{Scalarization of Set--Valued Functions}\label{sec:Scalarization}

\begin{definition}\label{def:scalarization}
  Let $g:X\to \P^\triup$ and $\phi:Z\to\OLR$ be two functions. The scalarization of $g$ with respect to $\phi$ is defined by
  \begin{align*}
    \forall x\in X:\quad \vp_{g,\phi}(x)=\inf\cb{-\phi(z)\,\st z\in g(x)}.
  \end{align*}
\end{definition}
If $\dom \phi=Z$, then  $\dom \vp_{g,\phi}=\dom g$.
Recall that in convex analysis, the indicator function of a set $M\subseteq X$ is the function $I_M:X\to\OLR$, $I_M(x)=0$, if $x\in M$ and $I_M(x)=+\infty$, else.
If $\phi\in Z^*$, then $\vp_{g,\phi}(x)=-\sigma(\phi|(g(x))$, the negative support function of $\phi$ at $g(x)$ and thus  $\vp_{g,0}(x)=I_{\dom g}(x)$.

If $g(x)=\cl\co(g(x))\in \G^\triup$, then by Equation \eqref{eq:separation}
\begin{align}\label{eq:descalarization}
\forall x\in X:\quad g(x)=\bigcap\limits_{z^*\in C^-\setminus\cb{0}}\cb{z\in Z\st \vp_{g,z^*}(x)\leq -z^*(z)}.
\end{align}
The scalarization $\vp_{g,0}$ can be omitted, as $\cb{z\in Z\st \vp_{g,0}(x)\leq 0}=Z$ holds for all $x\in \dom g$ and $\cb{z\in Z\st \vp_{g,z^*}(x)\leq -z^*(z)}=\emptyset$ for all $z^*\in C^-$ and $x\notin\dom g$.

A function $g:X\to\G^\triup$ is convex, positively homogeneous or subadditive, iff for all $z^*\in C^-$ the scalarization $\vp_{g,z^*}$ has the corresponding property. Closedness is not as immediate, as the following example shows. However, if all scalarizations $\vp_{g,z^*}$ with $z^*\in C^-\setminus\cb{0}$ are closed, then $g$ is closed.

\begin{example}\label{ScalarizationNotClosed}
  Let the set $Z=\R^2$ be ordered by the usual ordering cone $C=\R^2_+$, $z^*=(0,-1)$ and $g:\R\to\P^\triup$ be defined as $g(x)=\cb{(\frac{1}{x},0)}+C$, if $x>0$ and $g(x)=\emptyset$, else.
  Thus $\epi g$ is a closed set, while $\vp_{g,z^*}(0)=+\infty$ and $\vp_{g,z^*}(x)=0$ holds for all $x>0$ and therefore $\cl \vp_{g,z^*}(0)=0$.
\end{example}

The proof of the following Proposition follows \cite[Proposition 3.1]{EkelandTemam}.

\begin{proposition}\label{prop:OrderpreservationScalarization}
  Let  $g:X\to \P^\triup$ be a function, then
  \begin{align*}
        \forall x\in X:\quad \of{\cl \co g}(x)
        =\bigcap\limits_{z^*\in C^-\setminus\cb{0}}\cb{z\in Z\st \cl\co\vp_{g,z^*}(x)\leq -z^*(z)}.
  \end{align*}
\end{proposition}
\proof
For simplicity suppose that  $g$ is a closed convex function. The images $g(x)$ are elements of $\G^\triup$ and because of \eqref{eq:descalarization} it is left to prove
 \begin{align}\label{eq:DescalarizationClosedConvex}
   \forall x\in X:\quad g(x)
    \supseteq\bigcap\limits_{z^*\in C^-\setminus\cb{0}}\cb{z\in Z\st \cl\co\vp_{g,z^*}(x)\leq -z^*(z)}.
  \end{align}

If $g$ is improper, then $\dom g$ is closed and convex, thus $\dom g=\dom \cl\co\vp_{g,z^*}$ holds for all $z^*\in C^-$. If $z^*\in C^-\setminus\cb{0}$, then $\vp_{g,z^*}(x)=-\infty$ holds for $x\in \dom g$ and $\vp_{g,z^*}(x)=+\infty$, else. In this case, \eqref{eq:DescalarizationClosedConvex} is immediate.

Suppose $g$ is proper, $(x_0,z_0)\notin \epi g$. By Equation \eqref{eq:separation} it exists $(x^*,z^*)\in X^*\times Z^*$ and $t\in \R$ such that
\begin{align}\label{eq:proof3.3}
\forall (x,z)\in\epi g:\quad -x^*(x_0)-z^*(z_0)<t<-x^*(x)-z^*(z).
\end{align}

Thus $z^*\in C^-$ and $x^*_{-t}:X\to\R$ with $x^*_{-t}(x)=x^*(x)+t$ is an affine minorant of $\vp_{g,z^*}$, separating $(x_0,-z^*(z_0))$ from the epigraph of $\cl\co\vp_{g,z^*}$, proving
\begin{align*}
   \forall x\in X:\quad g(x)
    \supseteq \bigcap\limits_{z^*\in C^-}\cb{z\in Z\st \cl\co\vp_{g,z^*}(x)\leq -z^*(z)}.
  \end{align*}
The strict inequality in formula \eqref{eq:proof3.3} can be satisfied with $z^*=0$, if and only if $x_0\notin\dom g$. 

Next, chose $(\bar x,\bar z)\in X\times Z$ such that $\bar x\in \dom g$ and $\bar z\notin g(\bar x)$. As  $\bar x\in\dom g$, it exists $(\bar x^*,\bar z^*,\bar t)\in X^*\times (C^-\setminus\cb{0})\times \R$ with
\[
\forall (x,z)\in\epi g:\quad -\bar x^*(\bar x)-\bar z^*(\bar z)<\bar t< -\bar x^*(x)-\bar z^*(z).
\]
If $(\bar x^*,\bar z^*)$ separates $(x_0,z_0)$ from $\epi g$, then there is nothing more to prove. Otherwise, we can chose $s>0$ such that there exists $t_s\in\R$ with
\[
 -(x^*+s\bar x^*)(x_0)-(z^*+s\bar z^*)(z_0)<t_s< -(x^*+s\bar x^*)(x)-(z^*+s\bar z^*)(z)
\]
for all $(x,z)\in\epi g$. By assumption $(z^*+s\bar z^*)\in C^-\setminus\cb{0}$ is fulfilled and 
$(x^*+s\bar x^*)_{-t_s}$ separates $(x_0,-(z^*+s\bar z^*)(z_0))$ from $\cl\co\of{\epi \vp_{g,(z^*+s\bar z^*)}}$,
 thus $z_0\notin \cb{z\in Z\st \cl\co\vp_{g,(z^*+s\bar z^*)}(x_0)\leq -(z^*+s\bar z^*)(z)}$ and \eqref{eq:DescalarizationClosedConvex} is proven.
\pend

As an immediate corollary we get

\begin{corollary}\label{cor:ProperDescalarization}
  Let $f, g:X\to \P^\triup$  be two functions. Then  $(\cl\co f)\leq (\cl \co g)$ is satisfied, iff for all $z^*\in C^-\setminus\cb{0}$ the function $(\cl \co \vp_{f,z^*})$ is a minorant of $(\cl\co \vp_{g,z^*})$.

  Moreover, $\cl\co g:X\to\P^\triup$ is either proper or constant $Z$ if and only if the following equality is satisfied for all $x\in X$.
  \begin{align}\label{eq:ProperDescalarization}
        \of{\cl \co g}(x)=
        \bigcap\limits_{\substack{\cl\co \vp_{g,z^*}\text{ is proper},\\ z^*\in C^-\setminus\cb{0}}}
        \cb{z\in Z\st \cl\co\vp_{g,z^*}(x)\leq -z^*(z)}.
  \end{align}
\end{corollary}

Especially, $g:X\to\P^\triup$ is proper if and only if it is $z^*$--proper for some $z^*\neq 0$.

\begin{proposition}\label{prop:PropertiesOfScalarization}
  Let $I$ be an index set,  $f, g, g_i:X\to \P^\triup$ functions for all $i\in I$,  $T:X\to Y$ a linear continuous operator and $h:Y\to \P^\triup$ a function.
  Let $z^*\in C^-$, then the following formulas are true.
  \begin{enumerate}[(a)]
    \item
    $\forall x\in X:\quad \vp_{f+g, z^*}(x)=\vp_{f,z^*}(x)\isum \vp_{g,z^*}(x)$.

    \item
    $\forall x\in X:\quad \vp_{hT, z^*}(x)=\vp_{h,z^*}T(x)$.

    \item
    $\forall x\in X:\quad \vp_{\inf\limits_{i\in I} g_i, z^*}(x)=\inf\limits_{i\in I} \vp_{g_i,z^*}(x)$.

  \end{enumerate}
\end{proposition}
\proof
(a) and (b) are immediate from the Definition \ref{def:scalarization}. By definition $(\inf\limits_{i\in I} g_i)(x)=\bigcup\limits_{i\in I}g_i(x)$ 
for all $x\in X$. Thus, (c) is immediate.
\pend

Combining (a) and (c) from Proposition \ref{prop:PropertiesOfScalarization} one can derive
\begin{align}\label{eq:InfConvolOfScalarization}
\forall x\in X:\quad \vp_{f\isquare g, z^*}(x)=(\vp_{f,z^*}\isquare \vp_{g,z^*})(x)
\end{align}
and by (b) and (c)
\begin{align}\label{eq:Ag_OfScalarization}
\forall x\in X:\quad \vp_{Th, z^*}(x)=T\vp_{h,z^*}(x).
\end{align}

\begin{proposition}\cite{HamelSchrage10}\label{prop:InfResiduationOfScalarization}
  Let $f,g:X\to \P^\triup$ and $g_i;X\to \P^\triup$ be functions for all $i\in I$ and $x\in X$ and $z^*\in C^-$, then
  \begin{align*}
    \vp_{f,z^*}(x)\idif\vp_{g,z^*}(x) & \leq\vp_{\of{f\idif g},z^*}(x);\\
		\sup\limits_{i\in I} \vp_{g_i,z^*}(x) & \leq \vp_{\sup\limits_{i\in I} g_i, z^*}(x).
  \end{align*}
  If additionally  
 $f(x)=H_{\vp_{f,z^*}(x)}(z^*)$ holds true or $\sup\limits_{i\in I} g_i(x)=\bigcap\limits_{i\in I}H_{\vp_{g_i,z^*}(x)}(z^*)$ respectively,
 then we get equality.
\end{proposition}

If $\dual:(\P^\triup)^X\to (\P^\triup)^Y$ is a duality in the sense of Singer, \cite{Singer97} and 
$\dual(g)(x)=H_{\vp_{\dual(g),z^*}(x)}(z^*)$ holds true for all $g:X\to\P^\triup$
then by applying Proposition \ref{prop:PropertiesOfScalarization} (c) and Proposition \ref{prop:InfResiduationOfScalarization}  the following mapping is a duality:
\begin{align*}
 \forall g:X\to \P^\triup,\, \forall z^*\in C^-:\quad  \vp_{g,z^*}\mapsto \vp_{\dual(g),z^*}.
\end{align*}

\begin{definition}\label{def:Setification}
  Let $f:X\to\OLR$ and $\phi:Z\to\OLR$ be two functions. The set--valued function $S_{(f,\phi)}:X\to \P(Z)$ is defined by
  \[
  \forall x\in X:\quad S_{(f,\phi)}(x)=\cb{z\in Z\st f(x)\leq -\phi(z)}.
  \]
\end{definition}
A function $f_1$ is a minorant of $f_2:X\to\OLR$, iff $S_{(f_1,\phi)}(x)\supseteq S_{(f_2,\phi)}(x)$ is met for all $x\in X$ and all $\phi:Z\to\OLR$.

Each function $f:X\to\OLR$ is dominated by $\vp_{S_{(f,\phi)},\phi}$. If $\phi(X)\supseteq \R$, then equality holds true.

If $\phi:Z\to\OLR$ is nonincreasing, i.e. if $z_1\leq_C z_2$, then $\phi(z_2) \leq \phi(z_1)$, then $S_{(f,\phi)}$ takes its values in $\P^\triup$ for all $f:X\to\OLR$. If additionally $\phi=z^*\in C^-$, then the values of $S_{(f,z^*)}$ are  of the form 
$S_{(f,z^*)}(x)=H_{f(x)}(z^*)$, i.e. affine half--spaces or $\emptyset$ or $Z$, the image--space is a subset of $\P^\triup(Z,H(z^*))$. 

\begin{proposition}\label{prop:DominationProperty}
  Let $f:X\to\OLR$, $z^*\in C^-$ and $g:X\to\P^\triup$. Then  $S_{(f,z^*)}$ is a minorant of $g$ if and only if $f$ is a minorant of $\vp_{g,z^*}$.
\end{proposition}
\proof
First, let $S_{(f,z^*)}(x)\supseteq g(x)$ be assumed for all $x\in X$. Then
\[
\forall x\in X:\quad f(x)\leq \vp_{S_{(f,z^*)},z^*}(x)\leq \vp_{g,z^*}(x)
\]
is satisfied, $f$ is a minorant of $\vp_{g,z^*}$. On the other hand, if $f$ is dominated by $\vp_{g,z^*}$, then
\[
\forall x\in X:\quad S_{(f,z^*)}(x)\supseteq S_{(\vp_{g,z^*},z^*)}(x)\supseteq g(x).
\]
\pend

\begin{proposition}\label{prop:descalarizationConvexClosed2}
  Let $g:X\to \P^\triup$ be a function, $z^*\in C^-$. It holds
  \begin{align*}
  \forall x\in X:\quad S_{(\vp_{g,z^*},z^*)}(x) & =\cl\of{g(x)+H_0(z^*)};\\
  \cl \co (g(x)) & =\bigcap\limits_{z^*\in C^-\setminus\cb{0}}S_{(\vp_{g,z^*},z^*)}(x).
  \end{align*}
\end{proposition}
\proof
		By definition,
		\[
		\vp_{g,z^*}(x)=-\sigma(z^*|g(x)),\quad S_{(\vp_{g,z^*},z^*)}(x)=H_{0\isum (-\sigma(z^*|g(x))}(z^*),
		\]
		thus by Lemma \ref{lem:Review} the first equation is proven while the second is proven by a separation argument in $Z$, compare Equation \eqref{eq:separation}.
\pend

\begin{corollary}\label{cor:ImproperSetValuedConlinear}
  Let $f:X\to\OLR$ be a function $x\in X$. Then
  \begin{eqnarray*}
  f(x)\leq 0 \quad &\Leftrightarrow\quad S_{(f,0)}(x)=Z \quad &\Leftrightarrow\quad \vp_{S_{(f,0)},0}(x)=0;\\
  0 < f(x) \quad &\Leftrightarrow\quad S_{(f,0)}(x)=\emptyset \quad &\Leftrightarrow\quad \vp_{S_{(f,0)},0}(x)=+\infty.
  \end{eqnarray*}
 \end{corollary}
The scalarization $\vp_{S_{(f,0)},0}$ of $f:X\to\OLR$ is the indicator function of the sublevel set of $f$ at $0$.


\begin{example}\label{ex:ExtConlinearFunction}
  Let $x^*\in X^*$ and $r \in \R$. The function $x^*_r:X\to\R$ is defined by $x^*_r(x)=x^*(x)-r$ for all $x\in X$, the  closed improper $\inf$--extension $\hat x_r^*:X\to\R^\triup$ by
\[
\hat x_r ^*\of{x} = \left\{
       \begin{array}{lll}
         -\infty & : & x_r^*\of{x} \leq 0 \\
         +\infty & : & x_r^*\of{x} > 0
       \end{array}
     \right.
\]
If $x^*=0$, then we obtain  $\hat x^*_r=-\infty$ if $r\geq 0$, and $+\infty$ else.
If $f:X\to\R^\triup$ is a function, then
\begin{align}\label{eq:ImpropMinorantScalar}
  \forall x\in X:\quad \hat x^*_r(x)\leq f(x) \quad\Leftrightarrow\quad x^*(x)-r\leq I_{\dom f}(x).
\end{align}
We denote
\[
\forall z^*\in C^-,\; \forall r\in \R,\;\forall x\in X:\quad S_{(x^*,z^*,r)}(x)=S_{(x^*_r,z^*)}(x)
\]
and for completenes
\[
\forall z^*\in C^-,\;\forall x\in X:\quad S_{(x^*,z^*,+\infty)}(x)=Z,\qquad S_{(x^*,z^*,-\infty)}(x)=\emptyset.
\]
Thus, $S_{(x^*,z^*)}(x)=S_{(x^*,z^*,0)}(x)$ and for $(x^*,r)\in X^*\times \R$ and $z^*\in C^-$ it holds
\[
\forall x\in X:\quad S_{(\hat x^*_r,z^*)}(x)=S_{(x^*_r,0)}(x)=
    \left\{
      \begin{array}{ll}
        Z, & \hbox{if $x^*_r(x)\leq 0$;} \\
        \emptyset, & \hbox{else.}
      \end{array}
    \right.
\]
The scalarization $\vp_{S_{(x^*_r,0)},0}$ is the indicator function of $\dom \hat x^*_r$, while $\hat x^*_r=\vp_{S_{(\hat x^*_r,z^*)},z^*}$ is satisfied if  $z^*\neq 0$.
\end{example}

A function of the type $S_{(x^*,z^*,r)}:X\to \P^\triup$  with $(x^*,z^*,r)\in X^*\times C^-\times\R$ is called conaffine. If additionally $r=0$, then $S_{(x^*,z^*)}:X\to \P^\triup$ is called conlinear.
If $z^*\neq 0$, then $S_{(x^*,z^*,r)}$ is proper  and its values are affine half--spaces, 
if $z^*=0$, then $S_{(x^*,0,r)}(x)=Z$ for all $x\in\dom\hat x_r^*$ and $S_{(x^*,0,r)}(x)=\emptyset$ for $x\notin\dom \hat x^*_r$.

\begin{proposition}\label{prop:DualRepresentation}
  Let $g:X\to \P^\triup$, $z^*\in C^-\setminus\cb{0}$, $x^*\in X^*$ and $r\in \R$. Then $x^*_r$ is a minorant of $\vp_{g,z^*}$ if and only if $S_{(x^*,z^*,r)}$ is a minorant of $g$.
  The closed convex hull of $g$ is proper or constant $\emptyset$ or $Z$, if and only if
  \begin{align*}
  \forall x\in X:\quad (\cl \co g)(x)
    = \bigcap\limits_{\substack{x^*_r\leq \vp_{g,z^*},\\ z^*\in C^-\setminus\cb{0}}} S_{(x^*,z^*,r)}(x),
  \end{align*}
  it is improper, iff
  \begin{align*}
  \forall x\in X:\quad (\cl \co g)(x)
    = \bigcap\limits_{x^*_r\leq \vp_{g,0}} S_{(x^*,0,r)}(x).
  \end{align*}
\end{proposition}
\proof
By Proposition \ref{prop:DominationProperty} $x^*_r$ is a minorant of $\vp_{g,z^*}$ if and only if $S_{(x^*,z^*,r)}$ is a minorant of $g$
and by Corollary \ref{cor:ProperDescalarization} 
  \begin{align*}
  \forall x\in X:\quad
        \of{\cl \co g}(x) =\bigcap\limits_{z^*\in C^-\setminus\cb{0}}S_{(\cl\co\vp_{g,z^*},z^*)}(x)
  \end{align*}
  is satisfied. The closed convex hull of a scalar function $\vp_{g,z^*}:X\to\R^\triup$ is the supremum of its affine minorants if and only if $\cl\co\vp_{g,z^*}$ is proper or constant $+\infty$ or $-\infty$. 
 In this case,
  \[
  \forall x\in X:\quad S_{(\cl\co\vp_{g,z^*},z^*)}(x)= \bigcap\limits_{x^*_r\leq \vp_{g,z^*}}S_{(x^*_r,z^*)}(x).
  \]
  On the other hand, $\cl\co\vp_{g,z^*}$ is not proper for some $z^*\in C^-\setminus\cb{0}$, iff either $\vp_{g,z^*}\equiv +\infty$, thus $g\equiv \emptyset$, or $\vp_{g,z^*}$ does not have any affine minorants, proving the first statement.

  The function $\cl\co g:X\to\P^\triup$ is improper, iff $(\cl\co g)(x)=Z$ for all $x\in \cl\co\dom g=\dom (\cl\co g)$, iff for each $z^*\in C^-\setminus\cb{0}$ the scalarization $\vp_{g,z^*}$ is improper. It holds
  \[
  \forall x\in X:\quad \hat x^*_r(x)\leq \vp_{g,z^*}(x) \quad\Leftrightarrow\quad x^*_r(x)\leq I_{\dom \vp_{g,z^*}}(x)=\vp_{g,0}(x),
  \]
  compare \eqref{eq:ImpropMinorantScalar}. Thus
  \[
  \forall x\in X:\quad   x^*_r(x)(x)\leq \vp_{g,0}(x) \quad\Leftrightarrow\quad S_{(x^*,0,r)}(x)\supseteq g(x)
  \]
  and 
  \[
  \bigcap\limits_{x^*_r\leq \vp_{g,0}} S_{(x^*,0,r)}(x)=\begin{cases} 

  			Z,  & \text{ if } x\in\cl\co\dom g\\ 
  			\emptyset, & \text{else} 
																\end{cases} 
  \]
  proving the second statement.
\pend

 Especially,
  \[
  \forall x\in X:\quad  S_{(x^*,0,r)}(x)\supseteq g(x)
  \quad\Leftrightarrow\quad x^*(x)\idif r\leq I_{\dom g}(x).
  \]

Notice that the scalarizations of $g:X\to\P^\triup$ used in Proposition \ref{prop:DualRepresentation} are either proper or constant $+\infty$ and thus also the affine minorants needed in the representation are proper. Alternatively, the representation can be done excluding the $0$--scalarization $\vp_{g,0}$, compare \cite{HamelSchrage10} but at the cost of properness of the scalarizations and thus also of the affine minorants. The same is true for the following theorem, which is a corollary of Proposition \ref{prop:DualRepresentation} and Lemma \ref{lem:Separation}.

\begin{theorem}\label{thm:DualRepresentation}
Let $g:X\to \P^\triup$ be a function, then $g$ is convex and closed, if and only if $g$ is the point--wise supremum of its conaffine minorants. The function $g$ is proper or constant $\emptyset$ or $Z$, if and only if it is the point--wise supremum of its proper conaffine minorants,
\begin{align}\label{eq:dualRepresProper}
\forall x\in X:\quad g(x)
= \bigcap\limits_{\substack{S_{(x^*,z^*,r)}\lel_C g,\\ z^*\in C^-\setminus\cb{0}}} S_{(x^*,z^*,r)}(x).
\end{align}
Otherwise, $g$ is the point--wise supremum of its improper conaffine minorants,
\begin{align}\label{eq:dualRepresImproper}
\forall x\in X:\quad g(x)
= \bigcap\limits_{S_{(x^*,0,r)}\lel_C g} S_{(x^*,0,r)}(x).
\end{align}
\end{theorem}
Part of the statement of Theorem \ref{thm:DualRepresentation} has been proven for the proper case in \cite[Theorem 1]{Hamel09}, while in \cite[Theorem 5.30]{HamelSchrage10} a representation formula with improper conaffine minorants is proven.

\begin{example}\cite[Proposition 8]{Hamel09}
  Let $\bar g:X\to Z$ be a single--valued function with the epigraphical extension $g(x)=\cb{\bar g(x)}+C$ for all $x\in X$, $C$ a closed convex cone and $T:X\to Z$ is a linear continuous operator. Then
  \begin{align*}
    \forall x\in X:\quad \vp_{g,z^*}(x) & =-z^*(\bar g(x));\\
                   S_{(-T^*z^*,z^*)}(x) & =\cb{Tx}+ H(z^*)
  \end{align*}
  is satisfied for all $z^*\in C^-\setminus\cb{0}$. Moreover, $T(x)+z_0\leq  \bar g(x)$ is met for all $x\in X$ if and only if
    \begin{align*}
      \forall z^*\in C^-\setminus\cb{0},\; \forall x\in X:\quad
      S_{(-T^*z^*,z^*)}(x)+\cb{z_0}\supseteq g(x).
    \end{align*}
\end{example}

\begin{remark}
In \cite{KucukEtAl12JMAA}, the case of an ordering cone with compact base and nonempty interior is considered. The authors identify a vector--valued function $V_{f,K}:X\to Z$ with respect to a total ordering cone $K$. Under some additional assumptions,
\[
-z^*(V_{f,K}(x))=\vp_{f,z^*}(x)
\]
for all $z^*\in K^-$ and $V_{f,K}\in f(x)$. However, the existence of such a minimal element is only given under strong assumptions on both the ordering cone and the function itself.
\end{remark}

\section{Conjugation of Set--Valued Functions}\label{sec:Conjugation}

\begin{definition}\label{def:SetValConjugate}
The conjugate of a function $g:X\to \P^\triup$ is  $g^*:X^*\times C^-\times\R\to \P^\triup$, defined by
\begin{align*}
  \forall (x^*,z^*,r)\in X^*\times C^-\times\R:\quad
  g^*(x^*,z^*,r)=\sup\limits_{x\in X}\of{ S_{(x^*,z^*,r)}(x)\idif g(x) }.
\end{align*}
\end{definition}
For completenes we define
\begin{align*}
g^*(x^*,z^*,-\infty)&=\sup\limits_{x\in X}\of{ S_{(x^*,z^*,-\infty)}(x)\idif g(x) }=
		\begin{cases}  
			Z,\text{ if } \dom g=\emptyset,\\
			\emptyset,\text{ else,}
		\end{cases}\\
g^*(x^*,z^*,+\infty)&=\sup\limits_{x\in X}\of{ S_{(x^*,z^*,+\infty)}(x)\idif g(x) }=Z.
\end{align*}

 Especially, the values of $g^*$ are affine half--spaces or $\emptyset$ or $Z$.
With $(\vp_{g,z^*})^*:X^*\times\R\to\OLR$ we denote
\[
(\vp_{g,z^*})^*(x^*,r)=\sup\limits_{x\in X}\of{x^*_r(x)\idif\vp_{g,z^*}},
\]
compare \cite{HamelSchrage10} and abbreviate
$g^*(x^*,z^*)=g^*(x^*,z^*,0)$ and $(\vp_{g,z^*})^*(x^*)=(\vp_{g,z^*})^*(x^*,0)$ for all $(x^*,z^*)\in X^*\times C^-$.

\begin{proposition}\label{prop:ScalarizationOfConjugate}
  Let $g:X\to \P^\triup$ be a function and $x^*\in X^*$, $z^*\in C^-$ and $r\in\R$, then
 \begin{align*}
   \of{\vp_{g,z^*}}^*(x^*,r) &= \of{\vp_{g,z^*}}^*(x^*)\idif r,\\
   g^*(x^*,z^*,r)  &=H_{(\vp_{g,z^*})^*(x^*,r)}(z^*)
  \end{align*}
  is satisfied.
  If $z^*\in C^-\setminus\cb{0}$, then additionally
  \begin{align*}
  	\of{\vp_{g,z^*}}^*(x^*,r) &= \vp_{g^*(\cdot,z^*,\cdot),z^*}(x^*,r),\\
		g^*(x^*,z^*,r) &=g^*(x^*,z^*)\idif H_r(z^*)
  \end{align*}
	holds true while
  \begin{align*}\label{eq:0Conjugate}
        \of{\vp_{g,0}}^*(x^*,r)  &=\sigma(x^*|\dom g)\idif r,\\
        \vp_{g^*(\cdot,0,\cdot),0}(x^*,r)  &= I_{\dom g^*(\cdot,0,\cdot)}(x^*,r)
  \end{align*}
  and
  \begin{align*}
   g^*(x^*,0,r)  =
	\begin{cases} Z, &\text{if } \sigma(x^*,\dom g)\leq r\\
				\emptyset, &\text{else}
	\end{cases}
  \end{align*}
  is satisfied. 
\end{proposition}
\proof
The equation $\of{\vp_{g,z^*}}^*(x^*,r) = \of{\vp_{g,z^*}}^*(x^*)\idif r$ is proven by direct calculation, as by assumption $r\in\R$. It is well known and easy to prove that $(I_{\dom g})^*(x^*)=\sigma(x^*|\dom g)$, hence
\[
\of{\vp_{g,0}}^*(x^*,r)=(I_{\dom g})^*(x^*,r)=\sigma(x^*|\dom g)\idif r.
\]
Applying Lemma \ref{lem:Review} and Proposition \ref{prop:InfResiduationOfScalarization} gives
  \begin{align*}
   g^*(x^*,z^*,r)  =\bigcap\limits_{x\in X} H_{x_r^*(x)\idif\vp_{g,z^*}(x)}(z^*)= H_{(\vp_{g,z^*})^*(x^*,r)}(z^*)
  \end{align*}
and thus
	\begin{align*}
			\vp_{g^*(\cdot,z^*,\cdot),z^*}(x^*,r)  = \inf\cb{-z^*(z)\st (\vp_{g,z^*})^*(x^*,r)\leq -z^*(z) }.	\end{align*}
If additionally $z^*\neq 0$, then
	\begin{align*}
			\vp_{g^*(\cdot,z^*,\cdot),z^*}(x^*,r)  = \inf\cb{-z^*(z)\st (\vp_{g,z^*})^*(x^*,r)\leq -z^*(z) }= (\vp_{g,z^*})^*(x^*,r)
	\end{align*}
is satisfied while
  \begin{align*}
   g^*(x^*,0,r)  =H_{\sigma(x^*,\dom g)\idif r}(0)=
	\begin{cases} Z, &\text{if } \sigma(x^*,\dom g)\leq r\\
				\emptyset, &\text{else.}
	\end{cases}
  \end{align*}
and
	\begin{align*}
			\vp_{g^*(\cdot,0,\cdot),0}(x^*,r)  &
= \inf\cb{0\st (\vp_{g,z^*})^*(x^*,r)\leq 0 }\\
&=
	\begin{cases} 0, &\text{if } \sigma(x^*|\dom g)\idif r\leq 0\\
				+\infty, &\text{else}
	\end{cases}\\
&= I_{\dom g^*(\cdot,0,\cdot)}(x^*,r).
	\end{align*}

For all $z^*\in C^-\setminus\cb{0}$ and all $a,b\in \OLR$, by Lemma \ref{lem:Review} it holds
  \begin{align*}
   	H_{a\idif b}(z^*)= \cb{z\in Z\st a\leq -z^*(z)+b}=H_a(z^*)\idif H_b(z^*).
	\end{align*}
Thus especially $g^*(x^*,z^*,r)=g^*(x^*,z^*)\idif H_r(z^*)$ is true, proving the statement.  
\pend

Compare Corollary \ref{cor:ImproperSetValuedConlinear} for the scalarization with $z^*=0$. In Section \ref{sec:Duality}, set--valued duality results will be proven which only hold true on $X^*\times\cb{0}\subseteq X^*\times C^-$, if the pre--image space of the conjugate function is the set of conaffine, rather than conlinear functions. For this reason, we define the conjugate of a set--valued function to map from $X^*\times C^-\times\R$, rather than just $X^*\times C^-$ to $\P^\triup$.

\begin{proposition}\label{prop:ConjugateOfClosedConvexHull}
Let $g:X\to \P^\triup$ be a function, then the conjugate of $g$ and the conjugate of the closed convex hull of $g$ coincide,
    \[
    \forall x^*\in X^*, \forall z^*\in C^-, \forall r\in\R:\quad (\cl \co g)^*(x^*,z^*,r) = g^*(x^*,z^*,r).
    \]
\end{proposition}
\proof
By Proposition \ref{prop:OrderpreservationScalarization},
$\of{\cl\co \vp_{g,z^*}}(x)\leq \of{\cl\co \vp_{\cl\co g,z^*}}(x)$, thus
\begin{align}\label{eq:closedconvexScalarization}
\forall x\in X,\, \forall z^*\in C^-:\quad \cl\co \vp_{g,z^*}(x)= \cl\co \vp_{\cl\co g,z^*}(x)
\end{align}
and as
\begin{align*}
  \of{\cl\co g}^*(x^*,z^*,r)=\cb{z\in Z\st  \of{\vp_{\cl\co g,z^*}}^*(x^*)\leq r-z^*(z)}
\end{align*}
is satisfied for all $x^*\in X^*$, $r\in\R$ and all $z^*\in C^-$ and  we can conclude $\of{\cl\co g}^*(x^*,z^*,r)=g^*(x^*,z^*,r)$.
\pend

In contrast to the present approach, the (negative) conjugate in \cite{Hamel09} is defined as a $G^\triup$--valued function via an infimum rather than a supremum and thus avoiding a difference operation on the power set $\P\of{Z}$. In \cite{Diss} and in \cite{HamelSchrage10}, the same idea as in Definition \ref{def:SetValConjugate} has been used. In \cite{Diss}, the dual variables are reduced to the set $X^*\times C^-\setminus\cb{0}$, while in \cite{HamelSchrage10} the dual space is the set of all conaffine functions and again, $z^*=0$ is prohibited. There, improper scalarizations play an important role while in the present approach we avoid those at the expense of including $z^*=0$ in the dual space.

\begin{definition}\label{def:BiconjugateSetValued}
To a function $g:X\to \P^\triup$, the biconjugate $g^{**}:X\to \P^\triup$ is defined by
\begin{align}\label{eq:BiconjugateSetValued}
\forall x\in X:\quad g^{**}(x)=\bigcap\limits_{(x^*,z^*,r)\in X^*\times C^-\times \R}
    \of{ S_{(x^*,z^*,r)}(x)\idif g^*(x^*,z^*,r)}.
\end{align}
\end{definition}
It can be proven easily that
\[
\sup\limits_{(x^*,r)\in X^*\times \R}\of{ x^*_r(x)\idif (\vp_{g,z^*})^*(x^*,r) }=\sup\limits_{x^*\in X^*}\of{ x^*(x)\idif (\vp_{g,z^*})^*(x^*) }
\]
thus the biconjugate of a scalar function can be defined as usual, setting
\[
(\vp_{g,z^*})^{**}(x)=\sup\limits_{x^*\in X^*}\of{ x^*(x)\idif (\vp_{g,z^*})^*(x^*) }.
\]

\begin{theorem}[Biconjugation Theorem ]\label{thm:ScalarizationOfBiconjugate}
    Let $g:X\to \P^\triup$, then
  \begin{align}\label{eq:BiconjugateSetValuedProper}
  \forall x\in X:\quad (\cl\co g)(x)=g^{**}(x)=\bigcap\limits_{z^*\in  C^-}
				\cb{ z\in Z\st  (\vp_{g,z^*})^{**}(x)\leq -z^*(z)}
    \end{align}
    and $\cl\co g$ is proper or constant $\emptyset$ or $Z$, if and only if equality is satisfied when omitting $z^*=0$:
    \begin{align}\label{eq:ScalarizationBiconjProper}
      \forall x\in X:\quad & g^{**}(x)=\bigcap\limits_{z^*\in  C^-\setminus\cb{0}}
    				\cb{ z\in Z\st  (\vp_{g,z^*})^{**}(x)\leq -z^*(z)}.
    \end{align}
\end{theorem}
\proof
By Proposition \ref{prop:ScalarizationOfConjugate} the conjugate of a function is represented by $g^*(x^*,z^*,r)  =H_{(\vp_{g,z^*})^*(x^*)\idif r}(z^*)$ for all $(x^*,z^*,r)\in X^*\times C^-\times\R$, thus by Definition \ref{def:BiconjugateSetValued} we conclude
\[
 g^{**}(x)=\bigcap\limits_{z^*\in C^-}\of{ \bigcap\limits_{(x^*,r)\in X^*\times\R}
        \of{ H_{x^*(x)\idif r}(z^*)\idif H_{(\vp_{g,z^*})^*(x^*)\idif r}(z^*) }  }
\]
for all $x\in X$. Thus
\begin{align*}
 g^{**}(x)
&=\bigcap\limits_{z^*\in C^-} H_{\sup\limits_{x^*\in X^*}\of{x^*(x)\idif(\vp_{g,z^*})^*(x^*)}}(z^*)  \\
&=\bigcap\limits_{z^*\in C^-}\cb{z\in Z\st \vp^{**}_{g,z^*}(x)\leq -z^*(z)}
\end{align*}
is fulfilled for all $x\in X$, $\vp^{**}_{g,0}(x)=I_{\cl\co \dom g}(x)$ and
\[
 (\cl\co g)(x)=\bigcap\limits_{z^*\in C^-}\cb{ z\in Z\st  (\cl\co\vp_{g,z^*})(x)\leq -z^*(z)}.
\]
By the scalar biconjugation theorem, $(\cl\co\vp_{g,z^*})(x)=\vp^{**}_{g,z^*}(x)$ is met for all $x\in\cl\co\dom g$ and   $\cl\co\vp_{g,z^*}=\vp^{**}_{g,z^*}$ holds true iff $\cl \co \vp_{g,z^*}$ is either proper or constant $+\infty$ or $-\infty$, thus $(\cl\co g)(x)=g^{**}(x)$ for all $x\in X$.
If $\cl\co g$ is proper or constant $Z$ or $\emptyset$, then either $\cl\co\vp_{g,z^*}(x)=\vp_{g,z^*}^{**}(x)=+\infty$ for all $x\in X$ or by Corollary \ref{cor:ProperDescalarization}
  \begin{align*}
        \of{\cl \co g}(x)=
        \bigcap\limits_{\substack{\cl\co \vp_{g,z^*}\text{ is proper},\\ z^*\in C^-\setminus\cb{0}}}
        \cb{z\in Z\st \of{\vp_{g,z^*}}^{**}(x)\leq -z^*(z)}.
  \end{align*}
Moreover, $\of{\vp_{g,z^*}}^{**}$ is constant $-\infty$, whenever $\cl\co\of{\vp_{g,z^*}}$ is an improper function with $\dom\vp_{g,z^*}\neq\emptyset$, hence in this case
\[
\forall x\in X:\quad  g^{**}(x)=\bigcap\limits_{z^*\in  C^-\setminus\cb{0}}
    \cb{ z\in Z\st  (\vp_{g,z^*})^{**}(x)\leq -z^*(z)}
\]
is satisfied.
Finally, if \eqref{eq:ScalarizationBiconjProper} is met and $g^{**}(x)=Z$ for some $x\in X$, then
\[
\forall z^*\in C^-\setminus\cb{0}:\of{\vp_{g,z^*}}^{**}(x)=-\infty.
\]
In this case, $g^{**}$ is constant $Z$, proving the statement.
\pend

Assuming the order cone $C$ to be closed and pointed, a conjugate of a vector--valued function $f:X\to Z$ is defined in \cite{Borwein84a, Zalinescu83} and the references therein. The pre--image space of the conjugate is the set of continuous linear operators $T:X\to Z$, the conjugate is defined by
\[
f^+(T)=\sup\limits_{x\in X}\of{T(x)- f(x)}.
\]
To guarantee the existence of $f^+(T)$, the order induced by $C$ is assumed to fulfill a least upper bound property \cite{Zalinescu83} or even order completeness, \cite{Borwein84a}.
Identifying $f_C(x)= f(x)+C$, the following representation is fulfilled.
\begin{align*}
    f^+(T)+ H(z^*) & = \of{f_C}^*(-T^*z^*,z^*),\\
      f^+(T)+\cl C & = \bigcap\limits_{z^*\in C^-}\of{f_C}^*(-T^*z^*,z^*).
\end{align*}
Thus, results on the conjugate $f^+$ are included in the more general results on our set--valued conjugate.

The reader is referred to \cite[Proposition 12, 13; Theorem 2, 3]{Hamel09}, \cite[Section 4]{Diss}, \cite[Corollary 4.2]{GetanMaLeSi} for a more thorough investigation of the conjugates.

Theorem \ref{thm:ScalarizationOfBiconjugate} is a set--valued Fenchel--Moreau theorem, including the improper case alongside to the proper case.
The proper case can be found in \cite[Theorem 2]{Hamel09} or in \cite[Theorem 4.1.15]{Diss}.

\section{Duality Results}\label{sec:Duality}

In analogy to the scalar case, a chain--rule
as well as a Sandwich Theorem and the Fenchel--Rockafellar Duality Theorem  can be proven for set--valued functions.
We abbreviate the proofs by citing scalar results and applying
Proposition \ref{prop:ScalarizationOfConjugate} and Theorem \ref{thm:ScalarizationOfBiconjugate}.
Direct proofs for a special case can be found in \cite{Hamel09, Hamel10}.
There, strong duality results are formulated under the additional assumption of an inner point
$(x_0,z_0)\in\Int \epi g$ and $z^*\neq 0$.
We will show that continuity of $g$ in $x_0$ in the sense of \cite{AubinFrankowska90,Goepfert03},
too, is a sufficient assumption for strong duality results.

\begin{proposition}\label{prop:ContScalarization}
Let $z^*\in C^-$ and  $g:X\to \P^\triup$ a  convex function , $x_0\in\dom g$. If either of the following conditions is met, then  $\vp_{g,z^*}$ is convex and either continuous in $x_0$ or $\vp_{g,z^*}(x_0+x)=-\infty$ is satisfied for all elements $x$ of an open subset  $V\subseteq X$ with $0\in V$.

\begin{enumerate}[(a)]
	\item
	$g$ is  lower continuous in $x_0\in \dom g$ in the sense of  \cite[Definition 2.5.1.]{Goepfert03}, i.e. for all 	
	open sets $D\subseteq Z$ with $g(x_0)\cap D\neq\emptyset$ there exists  a $0$--neighborhood $V\subseteq X$
	such that
    \begin{align}\label{eq:ContinuousSetValuedFunctionLower}
        \forall x\in V:\quad g(x_0+x)\cap D\neq \emptyset;
    \end{align}

    \item
      there is $z_0\in g(x_0)$ such that $\cb{x\in X\st  z_0\in g(x)}$ is a neighborhood of $x_0$.
\end{enumerate}
\end{proposition}

\proof
As $g:X\to \P^\triup$  is by assumption convex, so is each scalarization $\vp_{g,z^*}$.
\begin{enumerate}[(a)]
  \item
To $z^*\in C^-$ and $t\in\OLR$, define the open set $S_t(z^*)=\cb{z\in Z\st  t<-z^*(z)}\subseteq Z$
and assume $g$ to be lower continuous in $x_0\in \dom g$ and $\varepsilon>0$. 

If $z^*=0$, then 
$g(x_0)\cap S_{-\varepsilon}(z^*)\neq \emptyset$
and thus it exists a $0$--neighborhood $V\subseteq X$ such that
$g(x_0+x)\cap S_{-\varepsilon}(z^*)\neq \emptyset$ is satisfied for all $x\in V$.
Therefore, $\vp_{g,0}=I_{\dom g}$ is  continuous at $x_0$.

If $z^*\in C^-\setminus\cb{0}$, then 
$g(x_0)\cap S_{\vp_{g,z^*}(x_0)+\varepsilon}(z^*)\neq \emptyset$ and it exists a $0$--neighborhood $V\subseteq X$ such that  $\vp_{g,z^*}$ is bounded from above on the set $\cb{x_0}+V$ by $(\vp_{g,z^*}(x_0)+\varepsilon)$ and thus  $\vp_{g,z^*}$ is either continuous in $x_0$ or $\vp_{g,z^*}(x_0+x)=-\infty$ for all $x\in V$.

\item
Let $g$ be convex and $z_0\in g(x_0)$ such that $N=\cb{x\in X\st  z_0\in g(x)}$ is a neighborhood of $x_0$, then
 $\vp_{g,z^*}$ is bounded from above on the set $N$ by $\vp_{g,z^*}(x)\leq -z^*(z_0)$ for all $x\in N$.
Thus  $\vp_{g,z^*}$ is either continuous in $x_0$ or $\vp_{g,z^*}(x_0+x)=-\infty$
 for all elements $x$ of an open set $V\subseteq X$ with $0\in V$.
\pend
\end{enumerate}

It is easy to check that under the assumptions of Proposition \ref{prop:ContScalarization} each scalarization satisfies $\vp_{g,z^*}(x_0)=\of{\cl\co\vp_{g,z^*}}(x_0)$. Hence the assumptions of Proposition \ref{prop:ContScalarization} are sufficient for the following two equalities
\begin{align}\label{eq:ResultOfContinuity}
g(x_0)=\of{\cl\co g}(x_0)=\cl\co\of{g(x_0)}
\end{align}

A more thorough investigation on continuity notions and sufficient constrained qualifications for strong duality in set--valued optimization will be done in the forthcoming work \cite{Heyde11W}. 

In the proof of Proposition \ref{prop:ContScalarization}, we use Property \eqref{eq:ContinuousSetValuedFunctionLower} for the open half--spaces $S_t(z^*)$ with $t\in\OLR$ and $z^*\in C^-$.
In case $Z$ is equipped with a norm, this leads to a uniform structure studied in \cite{SonntagZalinescu}, compare also \cite{LoehneZalinescu05, Loehne11Buch}.
Under the associated topology $\tau_S$ on $Z$, a sequence of closed convex sets $ \cb{A_n}_{n\in\N}\subseteq \G^\triup$ converges towards $A\in \G^\triup$ if and only if for all $z^*\in C^-$, $\liminf\limits_{n\to\infty}\cb{-z^*(a_n)\st  a_n\in A_n,}=\inf\cb{-z^*(a)\st  a\in A}$ is satisfied. Thus, if a convex function $g:X\to (\G^\triup,\tau_S)$ is continuous in $x_0\in X$ and $z^*$--proper, then especially  $\vp_{g,z^*}:X\to\R^\triup$ is continuous.
Also, if $g:X\to \G^\triup$ is continuous in the sense of \cite{Goepfert03} and $Z$ is equipped with a norm, then $g$ is continuous under the topology $\tau_S$.

\begin{definition}\label{def:InfConvForConjugate}
 Let $g_1, g_2:X\to \P^\triup$ be two functions,
$f:X\to \P^\triup$ a function, $z^*\in C^-$ and $T:X\to Y$, $S:Y\to X$ linear continuous operators.
\begin{enumerate}[(a)]
  \item
    Define the infimal convolution of $g^*_1$ and $g^*_2$ in $(x^*,z^*,r)\in X^*\times C^-\times\R$ with respect to $+$ by
    \begin{align*}
    (g_1^* \isquare g_2^*)(x^*,z^*)
	& =\cl\bigcup\limits_{\substack{x_1^*,x_2^*\in X^*,\,x_1^*+x_2^*=x^*\\r_1, r_2\in\R,\,r_1+r_2=r}}
    			\of{g^*_1(x_1^*,z^*,r_1)+ g^*_2(x_2^*,z^*,r_2)}.
    \end{align*}
  \item
    Define
    \begin{align*}
    \forall (x^*,r)\in X^*\times\R:\quad
    (T^*f^*)(x^*,z^*,r)&= \cl \bigcup\limits_{T^*y^*=x^*} f^*\of{y^*,z^*,r}\\
    (f^*S^*)(x^*,z^*,r)&= f^*\of{S^*x^*,z^*,r}.
    \end{align*}

\end{enumerate}
\end{definition}

As each image  $(g_1^* \isquare g_2^*)(x^*,z^*)$ and $(T^*f^*)(x^*,z^*,r)$ is by definition a closed half--space, we obtain  
\begin{align*}
(g_1^* \isquare g_2^*)(x^*,z^*)&=H_{-\sigma(z^*|(g_1^* \isquare g_2^*)(x^*,z^*))}(z^*);\\
(T^*f^*)(x^*,z^*,r)&=H_{-\sigma(z^*|(T^*f^*)(x^*,z^*,r))}(z^*).
\end{align*}
However, even with the images being closed it cannot be concluded that either function has a closed epigraph.

Applying the scalar chain--rule from Theorem \ref{thm:ScalarChainRule} and Propositions \ref{prop:PropertiesOfScalarization}, \ref{prop:InfResiduationOfScalarization}, \ref{prop:ScalarizationOfConjugate}, we get the following result.

\begin{theorem}[Chain--Rule ]\label{thm:ChainRule}
    Let $g:X\to\P^\triup$, $f:Y\to\P^\triup$ be two functions and $T:X\to Y$, $S:Y\to X$ linear continuous operators.

\begin{enumerate}[(a)]

\item
For all $(x^*,z^*,r)\in X^*\times C^-\times\R$ the conjugate of $\of{g\isquare Sf}$ is given by the function
\begin{align*}
\of{g\isquare Sf}^*(x^*,z^*,r)
&= g^*(x^*,z^*,r\idif\vp^*_{f,z^*}S^*(x^*))+ f^*S^*\of{x^*,z^*,\vp^*_{f,z^*}S^*(x^*)}\\
&\leq \inf\limits_{\substack{r_1+r_2=r,\\ r_1, r_2\in\R}} g^*\of{x^*,z^*,r_1}+f^*S^*\of{x^*,z^*,r_2}
\end{align*}
and equality holds true if $\dom g\neq\emptyset$ and $\dom f\neq\emptyset$.

\item
For all $(x^*,z^*,r)\in X^*\times C^-\times\R$ the conjugate of $\of{g+ fT}$ is is dominated as follows.
\begin{align*}
&\of{g+ fT}^*(x^*,z^*,r)\\
\leq &\inf\limits_{y^*\in Y^*} g^*(x^*-T^*y^*,z^*,r\idif\vp^*_{f,z^*}(y^*))+ f^*\of{y^*,z^*,\vp^*_{f,z^*}(y^*)}\\
\leq &\of{g^*\isquare T^*f^*}(x^*,z^*,r)
\end{align*}
and equality holds true in the second inequality if $\dom g\neq\emptyset$ and $\dom f\neq\emptyset$.

\item
    If $(fT)(x_0)+H(z^*)=Z$ for some $x_0\in \dom g$ or if both $f$ and $g$ are convex and one of the assumptions in Proposition \ref{prop:ContScalarization} is satisfied for $f$ in  an element of $T(\dom g)$, then for all $x^*\in X^*$ there exists an $y^*\in Y^*$ with 
\begin{align*}
       &\of{g+ fT}^*(x^*,z^*,r)\\ 
			=& \of{g^*\isquare T^*f^*}\of{x^*, z^*,r};\\
			=&  g^*\of{x^*-T^*y^*,z^*,r\idif\vp^*_{f,z^*}(y^*)}+ f^*\of{y^*,z^*,\vp^*_{f,z^*}(y^*)}.
\end{align*}

\end{enumerate}
\end{theorem}
\proof
\begin{enumerate}[(a)]
\item
By Proposition \ref{prop:ScalarizationOfConjugate}, the conjugate of a function $h:X\to\P^\triup$ can be represented as follows.
\[
\forall (x^*,z^*,r)\in X^*\times C^-\times\R:\quad h^*(x^*,z^*) = \cb{z\in Z\st  \vp^*_{h,z^*}(x^*)-r\leq -z^*(z)}.
\]
    Applying Propositions \ref{prop:PropertiesOfScalarization}, \ref{prop:InfResiduationOfScalarization} and the scalar chain--rule, Theorem \ref{thm:ScalarChainRule} we may conclude for $(x^*,z^*,r)\in X^*\times C^-\times\R$
\begin{align*}
\of{g\isquare Sf}^*(x^*,z^*,r)
         &= H_{\of{ \vp^*_{g,z^*}\ssum \vp_{f,z^*}^*S^*}\of{x^*}-r}(z^*)\\
	  &=\cb{z\in Z\st     \of{ \vp^*_{g,z^*}\ssum \vp_{f,z^*}^*S^*}\of{x^*}-r\leq -z^*(z) }.
 \end{align*}
If $\dom f=\emptyset$ or $\dom g=\emptyset$, then
\[
\of{ \vp^*_{g,z^*}\ssum \vp_{f,z^*}^*S^*}\of{x^*}-r=-\infty,
\]
and thus
\[
\of{g\isquare Sf}^*(x^*,z^*,r)
= g^*(x^*,z^*,r\idif\vp^*_fS^*(x^*))+ f^*S^*\of{x^*,z^*,\vp^*_fS^*(x^*)}=Z.
\]
Otherwise, both $\vp^*_{g,z^*}$ and $\vp_{f,z^*}^*S^*$ map into the set $\R\cup\cb{+\infty}$ and
\begin{align*}
&\of{ \vp^*_{g,z^*}\ssum \vp_{f,z^*}^*S^*}\of{x^*}-r\\
=&\of{ \vp^*_{g,z^*}\of{x^*}\idif \of{ r\idif \vp_{f,z^*}^*S^*\of{x^*}}}\isum \of{\vp_{f,z^*}^*S^*\of{x^*}\idif \vp_{f,z^*}^*S^*\of{x^*}}
\end{align*}
and a careful case study of $z^*=0$ and $\vp_{f,z^*}^*S^*\of{x^*}=+\infty$ gives
\begin{align*}
&\of{g\isquare Sf}^*(x^*,z^*,r)\\
		=& H_{ \vp^*_{g,z^*}\of{x^*}\idif \of{ r\idif \vp_{f,z^*}^*S^*\of{x^*}}}(z^*)+
		  	H_{\vp_{f,z^*}^*S^*\of{x^*}\idif \vp_{f,z^*}^*S^*\of{x^*}}(z^*)\\
		=&  g^*(x^*,z^*,\vp^*_{f,z^*}S^*(x^*))+ f^*S^*\of{x^*,z^*,\vp^*_{f,z^*}S^*(x^*)}.
\end{align*}
The second inclusion is immediate, if  $\dom f$ or $\dom g=\emptyset$. In case $\vp^*_{g,z^*}(x^*)=+\infty$ or $\vp^*_{f,z^*}S^*(x^*)=+\infty$, then for all $r_2\in\R$
\begin{align*}
	   & g^*(x^*,z^*,r\idif\vp^*_{f,z^*}S^*(x^*))+ f^*S^*\of{x^*,z^*,\vp^*_{f,z^*}S^*(x^*)}\\
	= & g^*(x^*,z^*,r-r_2)+ f^*S^*\of{x^*,z^*,r_2}=\emptyset.
\end{align*}
If  both $\vp^*_{g,z^*}(x^*)$ and $\vp^*_{f,z^*}S^*(x^*)\in\R$, then equality is proven by calculation.  

\item
    Applying Propositions \ref{prop:ScalarizationOfConjugate}, \ref{prop:PropertiesOfScalarization}, \ref{prop:InfResiduationOfScalarization} and the scalar chain--rule, Theorem \ref{thm:ScalarChainRule} we may conclude for $(x^*,z^*,r)\in X^*\times C^-\times\R$
\begin{align*}
(g+fT)^*(x^*,z^*,r)\supseteq H_{(\vp^*_{g,z^*}\ssquare T^*\vp^*_{f,z^*})(x^*)-r}
\supseteq H_{(\vp^*_{g,z^*}\isquare T^*\vp^*_{f,z^*})(x^*)-r}(z^*)
\end{align*}
and by the same arguments as above
\begin{align*}
&H_{(\vp^*_{g,z^*}\ssquare T^*\vp^*_{f,z^*})(x^*)-r}(z^*)\\
 \supseteq &\inf\limits_{y^*\in Y^*}\of{ g^*(x^*-T^*y^*,z^*,r\idif \vp_{f,z^*}^*(y^*)) +f^*(y^*,z^*,\vp_{f,z^*}^*(y^*))}\\
\supseteq &(g^*\isquare T^*f^*)(x^*,z^*,r) .
\end{align*}
Equality is shown by a case study with $ \vp_{f,z^*}^*(y^*)=+\infty$ when $\dom f\neq\emptyset$ and $\dom g\neq\emptyset$.

\item
Applying the scalar chain--rule,
\begin{align*}
(g+fT)^*(x^*,z^*,r)= H_{(\vp^*_{g,z^*}\isquare T^*\vp^*_{f,z^*})(x^*)-r}(z^*)= (g^*\isquare T^*f^*)(x^*,z^*,r)
\end{align*}
holds true under the given assumptions and for all $x^*\in X^*$ there is $y^*=T^*x^*$ such that
\begin{align*}
&(g+fT)^*(x^*,z^*,r)\\
= &g^*(x^*-T^*y^*,z^*,r\idif \vp_{f,z^*}^*(y^*))+ f^*(y^*,z^*,\vp_{f,z^*}^*(y^*)).
\end{align*}

\pend
\end{enumerate}

If $z^*\in C^-\setminus\cb{0}$, then 
\[
(g+fT)^*(x^*,z^*,r)\supseteq g^*\of{x^*-T^*y^*,z^*,0}+ f^*\of{y^*,z^*,0}+H_{-r}(z^*)
\]
and under the assumptions of Theorem \ref{thm:ChainRule} (c),
\[
(g+fT)^*(x^*,z^*,r)=g^*\of{x^*-T^*y^*,z^*,0}+ f^*\of{y^*,z^*,0}+H_{-r}(z^*)\neq Z
\]
 holds true. If additionally $\cl(fT(x_0)+H(z^*))=Z$, then $(g+fT)^*(x^*,z^*,r)=\emptyset$.

As in the scalar case, equality in Theorem \ref{thm:ChainRule} (a) and (b) does not hold true with the usual Minkowsky ($\inf$--) addition on the right hand side.
Indeed, if $g\equiv \emptyset$ and $fT(x_0)=Z$ for some $x_0\in X$, then for all $x\in X$ and all $(x^*,z^*,r) \in X^*\times C^-\setminus\cb{0}\times\R$ it holds
\begin{align*}
  &\of{g\isquare Sf}(x)=(g+fT)(x)=\emptyset;\\
  &g^*(x^*,z^*,r)=Z;\quad f^*S^*(x^*,z^*,r)=T^*f^*(x^*,z^*,r)=\emptyset.
\end{align*}
As $\emptyset$ dominates the Minkowsky sum, equality in general is not attained.
Notice however, that here as well as in the scalar case (see Theorem \ref{thm:ScalarChainRule} (d)) we do not assume properness for the strong chain--rule in Theorem \ref{thm:ChainRule} (c).

Setting $g=0$ or $X=Y$ and $S=T=id$, a sum--rule and a multiplication--rule are immediate corollaries of Theorem \ref{thm:ChainRule}.

\begin{corollary}[ Sandwich--Theorem ]
     Let  $T:X\to Y$ a linear continuous  operator  and $z^*\in C^-$. Let  $g:X\to \P^\triup$ and $f:Y\to \P^\triup$ be two convex functions 
such that
     \[
     \forall x\in X:\quad g(x)\subseteq H(z^*)\idif fT(x)
     \]
     and it exists $x_0\in \dom g$ such that one of the assumptions in Proposition \ref{prop:ContScalarization}
     applies for $f$ in $Tx_0\in Y$. Then there exists $y^*\in Y^*$ and $z_0\in Z$ such that
     \begin{align*}
     \forall x\in X:\quad & g(x)\subseteq S_{(T^*y^*,z^*)}(x)+\cb{-z_0}\subseteq H(z^*)\idif fT(x);\\
                          & z_0\in g^*(T^*y^*,z^*)\cap \of{ H(z^*)\idif f^*(-y^*,z^*) }.
     \end{align*}
     If additionally $\cl\of{g(x_0)+H(z^*)} = H(z^*)\idif fT(x_0)$ is fulfilled, then
     \[
     g^*(T^*y^*,z^*)= \cb{z_0}+ H(z^*);\quad  f^*(-y^*,z^*)= \cb{-z_0}+ H(z^*) .
     \]
\end{corollary}
\proof
By assumption, $\vp_{g,z^*},\vp_{f,z^*}:X\to\R^\triup$ are convex and proper,
$0\idif \vp_{f,z^*}\of{T(x)}\leq\vp_{g,z^*}(x)$
for all $x\in X$ and $\vp_{f,z^*}$ is continuous in $Tx_0$ and $x_0\in\dom \vp_{g,z^*}$. Thus,
\begin{align}\label{eq:SandwichProof1}
\forall x\in X:\quad 0\leq(\vp_{g,z^*}\isum\vp_{f,z^*}T)(x)
\end{align}
and it exists $y^*\in Y^*$ such that
\begin{align}\label{eq:SandwichProof2}
-\infty <s_0= (\vp_{g,z^*}+\vp_{f,z^*}T)^*(0)= \vp_{g,z^*}^*(T^*y^*)\isum\vp_{f,z^*}^*(-y^*).
\end{align}
By \eqref{eq:SandwichProof1}, $s_0\leq 0$ is valid and by \eqref{eq:SandwichProof2} it holds
\begin{align}\label{eq:SandwichProof3a}
\vp_{g,z^*}^*(T^*y^*)\leq -\vp_{f,z^*}^*(-y^*)\in\R
\end{align}
and thus the following holds true for all $x\in X$.
\begin{align*}
0\idif \vp_{f,z^*}(Tx)&\leq y^*(Tx)\idif \of{-\vp_{f,z^*}^*(-y^*)}\\
&\leq T^*y^*(x)\idif \vp^*_{g,z^*}(T^*y^*)\leq \vp_{g,z^*}(x),
\end{align*}
and if $\vp_{g,z^*}(x_0)=0\idif \vp_{f,z^*}(Tx_0)$, then $s_0=0$ and equality holds true.

Choose $z_0\in Z$ such that $-z^*(z_0)=-\vp_{f,z^*}^*(-y^*)\in\R$.
Applying Propositions \ref{prop:PropertiesOfScalarization}, \ref{prop:InfResiduationOfScalarization} and
\ref{prop:ScalarizationOfConjugate}, we get
\[
g^*(T^*y^*,z^*)\supseteq \cb{z_0}+H(z^*)=\of{ H(z^*)\idif f^*(-y^*,z^*)}
\]
and thus
\begin{align}\label{eq:SandwichProof4}
g(x)\subseteq S_{(T^*y^*,z^*)}(x)+\cb{-z_0}\subseteq H(z^*)\idif fT(x)
\end{align}

If $\cl\of{g(x_0)+H(z^*)} = H(z^*)\idif fT(x_0)$ holds, then $\vp_{g,z^*}(x_0)=\vp_{f,z^*}(Tx_0)$ is valid, thus
\[
 g^*(T^*y^*,z^*)=\cb{z_0}+H(z^*)
\]
and equality holds true in \eqref{eq:SandwichProof4}.
\pend

\begin{theorem}[Fenchel-Rockafellar-Duality ]\label{thm:RockafellarDuality}
  To $g:X\to \P^\triup$, $f:Y\to \P^\triup$ and a
  linear continuous operator $T:X\to Y$ and $z^*\in C^-$, denote
  \begin{align}
    P & = \cl \co \bigcup\limits_{x\in X}\of{ g(x)+f(Tx) };\\
    D(z^*) & = \bigcap\limits_{y^*\in Y^*}
                H(z^*)\idif \of{ g^*(T^*y^*,z^*)+ f^*(-y^*,z^*)}.
  \end{align}
\begin{enumerate}[(a)]
\item
   It holds $D(z^*)\supseteq P$.
\item
   If one of the assumptions in Proposition \ref{prop:ContScalarization}
   is in force for $f$ in an element in  $T(\dom g)$, then
    $\cl( P+H(z^*))=D(z^*)\neq \emptyset$ holds and it exists $y^*_{z^*}\in Y^*$ such that
  \[
  \cl(P+H(z^*)) = H(z^*)\idif \of{g^*(T^*y^*_{z^*},z^*)+ f^*(-y^*_{z^*},z^*)} \neq \emptyset.
  \]
In this case,  or if $fT(x_0)=Z$ for some $x_0\in \dom g$,
$P=\bigcap\limits_{z^*\in C^-\setminus\cb{0}}D(z^*)\neq \emptyset$ holds true and it exists a set
   $\cb{y^*_{z^*}\in Y^*\st  z^*\in C^-\setminus\cb{0}}$ such that
  \[
  P = \bigcap\limits_{z^*\in C^-\setminus\cb{0}} H(z^*)\idif \of{g^*(T^*y^*_{z^*},z^*)+ f^*(-y^*_{z^*},z^*)}.
  \]

\end{enumerate}
\end{theorem}

\proof
As $H(0^*)=Z$ and $D(0^*)=Z$, there is nothing left to prove for $z^*=0$.
If $z^*\in C^-\setminus\cb{0}$, then  the following inequality is met
\[
\sup\limits_{y^*\in Y^*}\of{0\idif \of{ \vp_{g,z^*}^*(T^*y^*)\isum \vp_{f,z^*}^*(-y^*)}}
\leq  \inf\limits_{x\in X}\of{\vp_{g,z^*}(x)+\vp_{f,z^*}T(x)}
\]
and equality holds, if $\vp_{f,z^*}$ and $\vp_{g,z^*}$ are proper functions and
$\vp_{f,z^*}$ is continuous in $Tx\in Y$ with $x\in \dom \vp_{g,z^*}$ or if either scalarization attains the value $-\infty$ within the domain of the other.
Applying Propositions \ref{prop:PropertiesOfScalarization}, \ref{prop:InfResiduationOfScalarization} and
\ref{prop:ScalarizationOfConjugate}  proves the statement.
\pend

We sum up our investigations by stating a set--valued version of the scalar Fundamental Duality Formula as can be found in \cite[Theorem 2.7.1]{Zalinescu02}.

\begin{theorem}[Fundamental Duality Formula ]\label{thm:FundDualityFormula}
Let $h:X\times Y\to \P^\triup$ be convex and such that there exists $x_0\in X$ with $(x_0,0)\in\dom h$. Let one of the
assumptions in Proposition \ref{prop:ContScalarization} be satisfied for $h(x_0,\cdot):Y\to \P^\triup$ in $0$.
\begin{enumerate}[(a)]
  \item
  If $h$ is $z^*$--proper for $z^*\in C^-\setminus\cb{0}$
, then
\begin{align*}
  \cl\bigcup\limits_{x\in X}\of{h(x,0)+ H(z^*)}=\bigcap\limits_{y^*\in Y^*}\of{H(z^*)\idif h^*(0,y^*,z^*)}
\end{align*}
and it exists $y_0^*\in Y^*$ such that 
\[
\cl\bigcup\limits_{x\in X}\of{h(x,0)+H(z^*)}=\of{H(z^*)\idif h^*(0,y_0^*,z^*)}.
\]
Furthermore,
\begin{align}\label{eq:z*-minimalelement}
\cl\of{h(\bar x,0)+H(z^*)}=\cl\bigcup\limits_{x\in X}\of{h(x,0)+H(z^*)}
\end{align}
holds for $\bar x\in X$ if and only if there is a  $\bar y^*\in Y^*$ such that
\begin{align}\label{eq:SubdiffFormula}
h^*(0,\bar y^*,z^*)\supseteq S_{((0,\bar y^*),z^*)}(\bar x,0)\idif h(\bar x,0).
\end{align}

  \item
  If $h$ is $z^*$--proper for all $z^*\in C^-\setminus\cb{0}$, then
\begin{align*}
  \cl\co\bigcup\limits_{x\in X}h(x,0)=\bigcap\limits_{\substack{y^*\in Y^*,\\ z^*\in C^-\setminus\cb{0}}}\of{H(z^*)\idif h^*(0,y^*,z^*)}
\end{align*}
and it exists a family $\cb{y^*_{z^*}\st  z^*\in C^-\setminus\cb{0}}\subseteq Y^*$ such that
\[
\cl\co\bigcup\limits_{x\in X}h(x,0)=\bigcap\limits_{z^*\in C^-\setminus\cb{0}}\of{H(z^*)\idif h^*(0,y^*_{z^*},z^*)}.
\]
 Furthermore,
\begin{align}\label{eq:minimalelement}
h(\bar x,0)=\cl\co\bigcup\limits_{x\in X}h(x,0)
\end{align}
holds for $\bar x\in X$ if and only if for all $z^*\in C^-\setminus\cb{0}$ there is a $\bar y^*_{z^*}\in Y^*$ such, that
\begin{align*}
\forall(x,y)\in X\times Y:\quad
h^*(0,\bar y^*_{z^*},z^*)\supseteq S_{((0,\bar y^*_{z^*}),z^*)}(\bar x,0)\idif h(\bar x,0).
\end{align*}

\end{enumerate}
\end{theorem}
\proof
\begin{enumerate}[(a)]
  \item
  If $z^*\in C^-\setminus\cb{0}$ and $h$ is $z^*$--proper or if $z^*=0$, we can derive that $\vp_{h,z^*}:X\times Y\to \R^\triup$
is convex and proper and
$\vp_{h,z^*}(x_0,\cdot):Y\to \R^\triup$ is continuous in $0$, compare Proposition \ref{prop:ContScalarization}.
Thus in both cases we can apply Theorem \ref{thm:ScalarFundamentalDuality} to  $\vp_{h,z^*}:X\times Y\to \R^\triup$ and attain
\begin{align*}
  \inf\limits_{x\in X}\of{\vp_{h,z^*}(x,0)}=\sup\limits_{y^*\in Y^*}\of{0\idif (\vp_{h,z^*})^*(0,y^*)},
\end{align*}
and the existence of $y^*\in Y^*$ such that $\inf\limits_{x\in X}\of{\vp_{h,z^*}(x,0)}=\of{0\idif (\vp_{h,z^*})^*(0,y^*)}$.
Furthermore,
\[
\vp_{h,z^*}(\bar x,0)=\inf\limits_{x\in X}\vp_{h,z^*}(x,0)
\]
is satisfied for $\bar x\in X$ if and only if $(0,\bar y^*)\in\partial \vp_{h,z^*}(\bar x,0)$, i.e.
\begin{align*}
\exists \bar y^*\in Y^*:\;  \forall(x,y)\in X\times Y:\quad
(\vp_{h,z^*})^*(0,\bar y^*)\leq (0,\bar y^*)(\bar x,0)\idif \vp_{h,z^*}(\bar x,0).
\end{align*}

To derive the set--valued result, we apply Propositions \ref{prop:PropertiesOfScalarization}, \ref{prop:InfResiduationOfScalarization}
and \ref{prop:ScalarizationOfConjugate} and achieve the desired.

  \item
  Assuming the $z^*$--properness for all $z^*\in C^-\setminus\cb{0}$, we obtain above results for all $z^*\in C^-\setminus\cb{0}$,
thus
\begin{align*}
  \cl\co\bigcup\limits_{x\in X}h(x,0) =  \bigcap\limits_{z^*\in C^-\setminus\cb{0}}\of{\cl\bigcup\limits_{x\in X}h(x,0)+ H(z^*)}
\end{align*}
and it exists $\cb{y^*_{z^*}}_{z^*\in C^-\setminus\cb{0}}\subseteq Y^*$ such that
\begin{align*}
\cl\co\bigcup\limits_{x\in X}\of{h(x,0)}
    =\bigcap\limits_{z^*\in C^-\setminus\cb{0}}\of{H(z^*)\idif h^*(0,y^*_{z^*},z^*)}.
\end{align*}

Finally, it exists a set
$\cb{y^*_{z^*}\st  z^*\in C^-\setminus\cb{0}}$ such that
\begin{align*}
\forall z^*\in C^-\setminus\cb{0}:\quad h^*(0,\bar y^*_{z^*},z^*)\supseteq S_{((0,\bar y^*_{z^*}),z^*)}(\bar x,0)\idif h(\bar x,0),
\end{align*}
if and only if
\[
\forall z^*\in C^-\setminus\cb{0}:\quad\cl\of{h(\bar x,0)+H(z^*)}=\cl\bigcup\limits_{x\in X}\of{h(x,0)+H(z^*)}
\]
or equivalently
\[
\forall z^*\in C^-\setminus\cb{0}:\quad\vp_{h,z^*}(\bar x,0)=\inf\limits_{x\in X}\vp_{h,z^*}(x,0)
\]
is satisfied. This is equivalent to
\begin{align*}
h(\bar x, 0) & = \bigcap\limits_{z^*\in C^-\setminus\cb{0}}\cl\of{h(\bar x,0)+H(z^*)}\\
             & = \bigcap\limits_{z^*\in C^-\setminus\cb{0}}\cl\bigcup\limits_{x\in X}\of{h(x,0)+H(z^*)}\\
             & \supseteq \cl\co\bigcup\limits_{x\in X}\bigcap\limits_{z^*\in C^-\setminus\cb{0}}\of{h(x,0)+H(z^*)}\\
             & = \cl\co\bigcup\limits_{x\in X} h(x,0)\supseteq h(\bar x,0).
\end{align*}
\pend
\end{enumerate}

The relation in Equation \eqref{eq:SubdiffFormula} is a set--valued variant of the subdifferential formula. 
In Equation \eqref{eq:minimalelement}, $h(\bar x,0)$  can be interpreted as a minimal point of the function $h(\cdot,0):X\to\P^\triup$, while in Equation \eqref{eq:z*-minimalelement} $\bar x$ is minimal with respect to the direction $z^*$.

\end{document}